\documentclass[10pt]{article}

\usepackage{amsfonts,latexsym,amstext}
\usepackage{amsmath}
\usepackage{amssymb}
\usepackage[english]{babel}
\usepackage[utf8]{inputenc}

\usepackage{hyperref}

\usepackage{enumerate}

\usepackage[usenames]{color}
\definecolor{red}{rgb}{1.0,0.0,0.0}

\definecolor{blu}{rgb}{0.0,0.0,1.0}

\definecolor{gre}{rgb}{0.03,0.50,0.03}

\definecolor{darkviolet}{rgb}{0.58, 0.0, 0.83}

\newtheorem{theorem}{Theorem}[section]
\newtheorem{lemma}[theorem]{Lemma}
\newtheorem{proposition}[theorem]{Proposition}
\newtheorem{definition}[theorem]{Definition}

\newtheorem{hypothesis}[theorem]{Hypothesis}

\newtheorem{remark}[theorem]{Remark}
\newtheorem{corollary}[theorem]{Corollary}
\numberwithin{equation}{section}
\def\qed{{\hfill\hbox{\enspace${ \square}$}} \smallskip}
\def\sqr#1#2{{\vcenter{\vbox{\hrule height .#2pt \hbox{\vrule
width .#2pt height#1pt \kern#1pt \vrule
width .#2pt} \hrule height .#2pt}}}}
\def\square{\mathchoice\sqr54\sqr54\sqr{4.1}3\sqr{3.5}3}

\def\qedo{\hbox{\hskip 6pt\vrule width6pt height7pt
depth1pt \hskip1pt}\bigskip}

\def\eps{\varepsilon}

\def\ds{\begin{displaystyle}}
\def\eds{\end{displaystyle}}
\def\dis{\displaystyle }
\def\<{\left\langle }
\def\>{\right\rangle }

\def\dim{\noindent \hbox{{\bf Proof.} }}
\def\R{\mathbb R}
\def\N{\mathbb N}

\def\C{\mathbb C}
\def\E{\mathbb E}
\def\P{\mathbb P}

\def\cala{{\cal A}}

\def\calc{{\cal C}}
\def\cald{{\cal D}}

\def\calf{{\cal F}}

\def\calh{{\cal H}}

\def\calk{{\cal K}}
\def\call{{\cal L}}
\def\caln{{\cal N}}

\def\calu{{\cal U}}

\def\call{{\cal L}}
\def\cals{{\cal S}}
\def\calo{{\cal O}}

\def\to{\rightarrow}

\begin{document}

\title{Lifting and partial smoothing for stationary HJB equations and related control problems in infinite dimensions}
\date{}

\author{Gabriele Bolli
\\
Dipartimento di Matematica Guido Castelnuovo,
Universit\`a degli Studi di Roma La Sapienza\\
Piazzale Aldo Moro 5,
00197 Roma,
Italy\\
e-mail: \texttt{gabriele.bolli@uniroma1.it}\\
\vspace{-0.6truecm}
\and
Fausto Gozzi\\
Dipartimento di Economia e Finanza, Universit\`a LUISS - Guido Carli\\
Viale Romania 32, 00197 Roma, Italy\\
e-mail: \texttt{fgozzi@luiss.it}
}

\maketitle

\vspace{-0.6truecm}

\begin{abstract}
{We study a family of stationary Hamilton-Jacobi-Bellman (HJB) equations in Hilbert spaces arising from stochastic optimal control problems. The main difficulties to treat such problems are: the lack of smoothing properties of the linear part of the HJB equation; the presence of unbounded control operators; the presence of state-dependent costs.
This features, combined together, prevent the use of the  classical mild solution theory of HJB equation (see e.g.
\cite[Ch.4]{FabbriGozziSwiech}).
The problem has been studied in the evolutionary case
in \cite{GozziMasiero23} using a "lifting technique" (i.e. working in a suitable space of trajectories where a "partial smoothing" property of the linear part of the HJB equations holds.
In this paper we extend such a theory to the case of infinite horizon optimal control problems, which are very common, in particular in economic applications.
The main results are: the existence and uniqueness of a regular mild solution to the HJB equation; a verification theorem, and the synthesis of optimal feedback controls.}
\end{abstract}


\smallskip \noindent


\smallskip \noindent

\textbf{Acknowledgements}:
We thank Federica Masiero for useful discussions. We also thanks the seminar participants in Durham and Wien for useful questions and remarks.

\tableofcontents

\section{Introduction}
{This paper deals with the solution of a family of stationary Hamilton-Jacobi-Bellman (HJB) equations arising from infinite horizon stochastic optimal control problems in infinite dimensions and their applications. A standard method for solving such equations is to formulate them in an integral form (the so-called mild form) and use a fixed-point argument, a survey on this can be found, e.g., in
\cite[Chapter 4]{FabbriGozziSwiech}). This approach, however, crucially relies on the regularizing properties of the transition semigroup associated with the underlying stochastic process.}

{Hence when such a smoothing property is missing, as it happens in many applications (like the ones involving, as state equations, delay equations or age structured PDEs), the problem becomes much more challenging from the technical viewpoint. Things get worse when this feature is coupled with the presence of unbounded control operators and/or with the presence of state dependent objective function. This kind of problem have been studied, only for finite horizon control problems, in cases of increasing difficulties, in the papers
\cite{FGFM-I,FGFM-II,FGFM-III,GozziMasiero23}; in particular the last one considers the case where all the above three features arise together.
The core part of this last paper is to lift the state space into a suitable space of trajectories (what we call here
"lifting technique") and to show that, in such spaces, a suitable
"partial smoothing" property holds.}

%running cost of the problem depends on the state variable, especially when the associated Ornstein-Uhlenbeck semigroup lacks adequate regularizing properties. In this scenario, which is common in applications, the fixed-point strategy often fails.
%This difficulty is particularly acute in the infinite-horizon setting, where integration over time does not provide the smoothing effect sometimes available in finite-horizon problems. Our work aims to overcome these issues, also allowing for the presence of unbounded control operators.

\medskip

{The central contribution of this work is the extension of the above methodology ("lifting technique" and "partial smoothing") to the infinite-horizon setting:
\vspace{-0.2truecm}
\begin{itemize}
  \item
finding solutions of the HJB equations that are regular enough for the candidate optimal feedback map to be well-defined.
\vspace{-0.2truecm}
  \item
using such result to prove verification theorems and existence of optimal feedback controls.
\vspace{-0.2truecm}
  \end{itemize}
Due to the complexity of such a methodology, such extension is delicate and nontrivial, see Remarks \ref{rem:space}, \ref{rem:fixedpoint}, and for details.}

\medskip

{Our main motivation is that problems displaying the above features are frequently encountered in the modeling of applied systems, like the cases of problems with boundary control. Our primary examples here will be a stochastic wave equation with distributed control and a stochastic heat equation with boundary control. We are currently working on extensions to cover problems with delay in the control and problem with age structured state equations, but this requires additional results which will be the object of a future paper.}

%This method, recently developed for finite-horizon problems in \cite{GozziMasiero23}, is specifically designed to overcome the lack of regularizing properties of the semigroup. The core idea is to shift the problem from a space of functions of the state, $\phi(x)$, to a space of functions of the entire trajectory generated by the state, $\hat{\phi}(y_x(\cdot))$. Since trajectories often exhibit better regularity properties than single points, this "lifted" setting can restore the necessary smoothing properties, making the fixed-point argument viable.

{We must add that these kind of HJB equations could be in principle studied by other methods, like the viscosity solution approach (see a survey on it in \cite[Ch.3]{FabbriGozziSwiech})
or the backward SDE approach (see a survey on it in \cite[Ch.6]{FabbriGozziSwiech}). However, up to now, the viscosity solution approach does not provide regularity results
in the case we are interested in, in particular in presence of unbounded control operator; on the other hand
the backward SDE approach relies on structural conditions of the problem which, again, do not hold in most of the applied examples we have in mind, in particular when the control operator
is unbounded (like in boundary control problems).}

\medskip

The structure of the paper is as follows:
\begin{itemize}
\vspace{-0.2truecm}
\item In Section \ref{sec:abstract_formulation}, we set up the abstract framework for the infinite-horizon stochastic optimal control problem, accommodating also unbounded control operators.
\vspace{-0.2truecm}
\item In Section \ref{sec:smoothing_and_lifting}, we recall the partial smoothing results of the papers \cite{FGFM-I,FGFM-III} and adapt them to our infinite horizon case. Moreover we also recall and adapt the "lifting" technique of \cite{GozziMasiero23}.
\vspace{-0.2truecm}
\item In Section \ref{sec:hjb_equation}, we use the lifting framework to study the stationary HJB equation. We prove the existence and uniqueness of a regular mild solution via a fixed-point argument for a sufficiently large discount factor.
\vspace{-0.2truecm}
\item In Section \ref{sec:verification_and_synthesis}, we establish a verification theorem that identifies the HJB solution with the value function and allows for the synthesis of an optimal feedback control.
\vspace{-0.2truecm}
\item In Section \ref{sec:wave_application}, we apply our theory to a stochastic wave equation with distributed control.
\vspace{-0.2truecm}
\item In Section \ref{sec:application}, we apply our theory to the stochastic heat equation with Dirichlet boundary control.
\end{itemize}

\section{Abstract Formulation of our Control Problem}\label{sec:abstract_formulation}
We now introduce the abstract formulation of the family of infinite horizon stochastic optimal control problems treated in this paper.
The abstract framework is the one of stochastic optimal control in infinite dimension 
(see e.g. \cite[Chapter 2]{FabbriGozziSwiech}) with unbounded control operators which will include 
various problems interesting for applications: in particular Boundary Control of Stochastic Partial Differential Equations (SPDEs) and controlled Stochastic Delay Differential Equations (SDDEs) with delay in the control.
\\
We start with the following assumptions about the spaces and the probabilistic framework.
\begin{hypothesis}\label{ip:spaces}
\begin{itemize}
%\item []
\item [(i)] The state space $H$ the control space $K$ and the noise space $\Xi$ are
real separable Hilbert spaces.
\item[(ii)] $\overline{H}$ (the space containing the image of the control operator) is a real separable Banach space such that $H\subseteq \overline{H}$ with continuous and dense inclusion.
\item[(iii)]
$\left(\Omega, \calf, \P\right)$ is a complete probability space.
\end{itemize}
\end{hypothesis}
% $H$ and a real separable Banach space $\overline{H}$ such that $H\subseteq \overline{H}$ with continuous and dense inclusion.
%We identify $H$ with its dual and we get a so-called Gelfand triple
%\begin{equation}
%V'\hookrightarrow H \hookrightarrow V.
%\end{equation}
\begin{remark}
%This Gelfand triple structure is a standard setting for the analysis of SPDEs, particularly when dealing with unbounded operators \cite{{FabbriGozziSwiech}}, \cite{{DP14}}.
In the boundary control case mentioned in the introduction the state space $H$ is $L^2(\mathcal{O})$ while the larger Banach space $\overline{H}$ will be a negative Sobolev space $H^{-\alpha}(\mathcal{O})$ with $\alpha$ properly chosen (see \cite[Appendix C]{FabbriGozziSwiech})
is introduced since the, possibly unbounded, control operator $B$ takes values in $\overline{H}$ and not necessarily in $H$. This is exactly what happens in the case of boundary control problems and in the case of pointwise delay in the control, see e.g. \cite[Section 3]{FGFM-III}, \cite[Section 4]{GozziMasiero23}.
\hfill\qedo
\end{remark}

\vspace{-0.4truecm}

The state equation is the following evolution equation on $\overline{H}$:
\begin{equation}\label{equazione stato}
\left\{
\begin{array}{l}
\dis
d X(s)= AX(s)\,ds+ Bu(s)\,ds +G\,dW(s), \quad s\in(0,\infty)
\\\dis
X(0)=x\in H.
\end{array}
\right.
\end{equation}
%This equation represents the abstract form of a controlled linear stochastic evolution system. The main mathematical challenge that this paper aims to address arises from the properties of the control operator $B$, which, in many significant applications, is unbounded.

Here we make the following assumptions.
\begin{hypothesis}\label{ipo-data}
\begin{itemize}
%\item []
\item [(i)] $A$ is the infinitesimal generator of a strongly continuous semigroup $\left\{e^{tA}\right\}_{t\geq 0}$ on $H$. We assume that such a semigroup can be extended to a strongly continuous semigroup on $\overline{H}$ and denote this extension by $\left\{\overline{e^{tA}}\right\}_{t\geq 0}$ \footnote{This property is satisfied in the boundary control of the stochastic heat equation, as explained in detail in section \ref{sec:application}.}
\vspace{-0.2truecm}
\item [(ii)] $B\in \call(K, \overline{H})$.
\item [(iii)] $G\in\call(\Xi,H)$.
\vspace{-0.2truecm}
\item [(iv)]
 $W$ is an
$\left(\Omega, \calf, \left(\calf_t\right)_{t\geq 0}, \P\right)$-cylindrical Wiener process in $\Xi$, and $\left(\calf_t\right)_{t\geq 0}$ is the augmented filtration generated by $W$.
\item [(v)]  
The selfadjoint operator
$Q_t:=\int_0^t e^{sA}GG^*e^{sA^*}\,ds$ is trace class for all $t>0$.
\item [(vi)]
The space of admissible control is given by
(here $U$ is a closed and bounded subset of $K$)
 \vspace{-0.2truecm}
\begin{equation*}\label{eq:admcontr}
\calu:=\left\{
u:[0,\infty)\times \Omega \to U \subseteq K, \left(\calf_t\right)_{t\geq 0} \text{ progressively measurable}
\right\}.
\vspace{-0.2truecm}
\end{equation*}
\end{itemize}
\end{hypothesis}
\begin{remark}
Hypothesis $(i)$ is the standard assumption for modeling linear evolution processes over time {\cite[Appendix A]{DaPratoZabczyk14}}. Hypothesis $(ii)$ formalizes the unboundedness of the control operator. The operator $B$ maps the control space $K$ not into the state space $H$, but into the larger space $\overline{H}$. This setup is essential for problems like boundary control, where the control's effect is too "rough" to remain in $H$ \cite{FGFM-III}, while in Hypothesis $(v)$ the trace-class condition on the covariance operator $Q_t$ of the noise is a fundamental requirement to ensure that the stochastic convolution term in the mild solution is a well-defined Gaussian process with values in $H$ {\cite[Chapter 5, Section 5.1.2]{DaPratoZabczyk14}}.
\hfill\qedo
\end{remark}

\vspace{-0.4truecm}

Equation (\ref{equazione stato}) is formal and has to be considered in its mild formulation.
\begin{definition}
We say that $X(\cdot)$ is a mild solution of (\ref{equazione stato}) if
\begin{equation}\label{mild}
X(s)=\overline{e^{sA}}x+\int_0^s\overline{e^{(s-r)A}}B u(r)\,dr +\int_0^se^{(s-r)A}G\,dW(r),
\text{ \ \ \ }s\in[0,\infty).
\end{equation}
We denote by $X(t; x,u)$ the mild solution of $(\ref{equazione stato})$ at time $t\geq0$ with initial condition $x\in H$ and control $u\in \cal U$.
\end{definition}
\begin{remark}
The mild solution is an integral representation of the state equation. This formulation is essential because classical (strong) solutions may not exist due to the irregularity of the noise and control terms. The term $\int_0^s\overline{e^{(s-r)A}}B u(r)\,dr$ is particularly critical; it requires extending the semigroup $e^{tA}$ to the larger space $\overline{H}$ to properly handle the action of the unbounded control operator $B$.
\end{remark}
Our goal is to minimize, over the set of admissible control strategies $\mathcal{U}$, the following classical infinite-horizon discounted cost functional, defined for $x\in H$
\vspace{-0.2truecm}
\begin{equation}\label{costo astratto}
J(x;u)=\E \left(\int_0^\infty e^{-\lambda s}\left[\ell_0\left(X(s; x, u)\right)+\ell_1\left(u(s)\right)\right]\,ds \right).
\vspace{-0.2truecm}
\end{equation}
where, as usual, the factor $e^{-\lambda s}$ with $\lambda > 0$ ensures the convergence of the integral and gives more weight to costs incurred in the near future.
We make the following assumptions on the cost functions $\ell_0$ and $\ell_1$.
\begin{hypothesis}
\begin{itemize}
\item [(i)] $\ell_0: \overline{H}\rightarrow \R$ is measurable, bounded from below and satisfies the growth condition $\left|\ell_0(x)\right|\leq C_{0}\left(1+\left|x\right|^{p}\right)$,
for some positive costants $C_{0}$, $p$ and any $x\in \overline{H}$.
\item [(ii)] $\ell_1: U\rightarrow \R$ is measurable and bounded.
\end{itemize}
\end{hypothesis}
Under this assumptions, the functional $J$ is well defined and bounded from below for all $x\in \overline{H}$, so we define the value function of the problem as
\vspace{-0.2truecm}
\begin{equation}\label{funzione valore}
V(x):= \inf_{u\in \calu}J(x;u),
\qquad x \in \overline{H}.
\vspace{-0.2truecm}
\end{equation}
Our goal here is to follow the dynamic programming approach whose core is to characterize this value function $V$ as the unique solution to the associated Hamilton-Jacobi-Bellman (HJB) equation, see e.g. \cite[Chapter 2]{FabbriGozziSwiech}. 
\begin{remark}
We choose the current cost $\ell_{0}$ to be defined on the extended space $\overline{H}$. In this way, the functional $J(x;u)$ is well defined also for $x\in \overline{H}$ and the value function can be studied directly on $\overline{H}$.
This is a subtle but important modeling choice. Defining the cost on the larger space $\overline{H}$ is consistent with the possibility that the state trajectory $X(s)$ may not remain in $H$, which is a key feature of systems with unbounded control operators. 
However, if for any control $u\in \cal U$ the integral term $\int_0^t\overline{e^{(t-s)A}}B u(s)\,ds$ belongs to $H$ for all $t\geq0$, this implies that the solutions of (\ref{equazione stato}) belong to $H$ for all times whenever the initial condition does. Hence, under this additional assumption, one can see the value function defined on (\ref{funzione valore}) as an extension of a well defined value function, defined in the same way as (\ref{funzione valore}), but in the smaller space $H$.
\end{remark}

\section{Recalling and adapting Partial Smoothing results}\label{sec:smoothing_and_lifting}

\subsection{''Basic'' partial smoothing}
\label{sec:partial_smoothing}

%This section is dedicated to the study of the regularizing properties of the Ornstein-Uhlenbeck (O-U) semigroup, which is the transition semigroup associated with the uncontrolled part of our state equation. The properties of this semigroup are fundamental for solving the Hamilton-Jacobi-Bellman (HJB) equation associated with the control problem. Standard results for HJB equations typically require the semigroup to be strongly Feller, a property that is often too restrictive for various systems \cite[Section 3]{FGFM-III}.
%\\
%a weaker but still useful property known as
We now recall and adapt the results of "partial smoothing" given in In \cite{FGFM-I}. 
%The core idea is that, although the semigroup may not regularize functions in all directions, it may do so for a certain class of functions in the specific directions related to the control action and this may be enough to solve our HJB equations.
%\\
%Then, in what follows, we recall the results of \cite{FGFM-I} for the commodity of the reader.
We introduce the Ornstein-Uhlenbeck process on $H$ as the mild solution of the following uncontrolled equation on $H$
\begin{equation}\label{ou-pro}
\left\lbrace\begin{array}{l}
dZ(t)=AZ(t)\,dt+G\,dW(t), \quad t\geq0\\
X(0)=x\in H.
\end{array}\right.
\end{equation}
and we denote it by $X(\cdot;x)$. The associated Ornstein-Uhlenbeck semigroup is defined on $B_{b}(H)$ by:
\begin{equation}\label{ou-sgr}
P_{t}[\phi](x):= \E\left[\phi\left(X(t;x)\right)\right]= \int_{H}\phi\left(e^{tA}x+y\right)\,\caln(0,Q_{t})(dy), \quad x\in H, \quad t\geq 0.
\end{equation}

In the same way, if we denote by $\overline{A}$ the infinitesimal generator of the extended semigroup $\overline{e^{tA}}$ then we can define the extended Ornstein-Uhlenbeck process in the extended state space $\overline{H}$ as the mild solution of
\begin{equation}\label{ou-pro-ext}
\left\lbrace\begin{array}{l}
dX(t)=\overline{A}X(t)\,dt+G\,dW(t), \quad t\geq0\\
Z(0)=z\in \overline{H}.
\end{array}\right.
\end{equation}
The associated semigroup on $B_{b}(\overline{H})$ is defined similarly as
\begin{equation}\label{ou-sgr-ext}
\overline{P}_{t}[\xi](z):= \E\left[\xi\left(Z(t;z)\right)\right]= \int_{\overline{H}}\xi\left(\overline{e^{tA}}z+y\right)\,\caln(0,\overline{Q_{t}})(dy), \quad z\in \overline{H}, \quad t\geq 0,
\end{equation}
where $\caln(0,\overline{Q}_{t})$ is the law of the gaussian process $W^{\overline{A}}(t)=\int_0^t \overline{e^{(s-r)A}}G\,dW(r)$ defined on $\overline{H}$ whose covariance operator $\overline{Q}_{t}\in \call(\overline{H},\overline{H}')$ is given by $\overline{Q}_{t}=\int_{0}^{t}e^{s\overline{A}}GG^{*}e^{s\overline{A^{*}}}\,ds$.\newline

The following computation shows that the semigroup $\overline{P}_{t}$ can be seen as an extension on $B_{b}(\overline{H})$ of the semigroup $P_{t}$ naturally defined on $B_{b}(H)$.
\begin{lemma}\label{lemma-ou}
Let $\phi\in B_{b}(\overline{H})$. Let $\tilde{\phi}$ be the restriction of $\phi$ to $H$, which is an element of $B_{b}(H)$ thanks to the continuous inclusion $H \hookrightarrow \overline{H}$. Let $\{P_{t}\}_{t\geq0}$ and $\{\overline{P}_{t}\}_{t\geq0}$ be the Ornstein-Uhlenbeck semigroups defined on (\ref{ou-sgr}) and (\ref{ou-sgr-ext}) respectively. Then for each $x\in H$ and $t\geq 0$ one has
\begin{equation*}
P_{t}[\tilde{\phi}](x)=\overline{P}_{t}[\phi](x).
\end{equation*}
\end{lemma}

\dim
By definition (\ref{ou-sgr}) we have $\overline{P}_{t}[\phi](x):= \E\left[\phi\left(Z(t;x)\right)\right]$, being $Z(\cdot;x)$ the solution of (\ref{ou-pro-ext}) with initial datum $x\in H$. Moreover, being $\overline{e^{tA}}$ an extension of $e^{tA}$ it follows immediatly that $Z(t;x)=X(t;x)$ for all $t\geq0$, where $X(\cdot;x)$ is the solution of (\ref{ou-pro}) with the same initial condition. Since $\tilde{\phi}$ is the restriction of $\phi$ to $H$ it is clear that $\phi\left(X(t;x)\right)=\tilde\phi\left(X(t,x)\right)$ and hence $\E\left[\phi\left(Z(t;x)\right)\right]= \E\left[\tilde{\phi}\left(X(t;x)\right)\right]$. By definition the last term coincide with $P_{t}[\tilde{\phi}](x)$.
\hfill\qed

As a consequence of the previous lemma we have the following representation of the extended Ornstein-Uhlenbeck semigroup.
\begin{corollary}
Let $\phi\in B_{b}(\overline{H})$. Let $\{\overline{P}_{t}\}_{t\geq 0}$ be the extended Ornstein-Uhlenbeck semigroups defined on (\ref{ou-sgr-ext}). Then for each $x\in \overline{H}$ and $t\geq 0$ the following equality holds
\begin{equation*}
\overline{P}_{t}[\phi](x)= \int_{H}\phi\left(\overline{e^{tA}}x+y\right)\,\caln(0,Q_{t})(dy).
\end{equation*}
\end{corollary}
\dim
Thanks to (Lemma \ref{lemma-ou}) we have that for all $x\in H$
\begin{equation}\label{eq1}
\overline{P}_{t}[\phi](x)= P_{t}[\tilde\phi](x)= \int_{H} \tilde\phi\left(e^{tA}x+y\right)\,\caln(0,Q_{t})(dy)
\end{equation}
being $\tilde\phi$ the restriction to $H$ of $\phi$. Clearly, the last term of equation (\ref{eq1}) is the restriction to $H$ of the continuous function $g(x)=\int_{H} \phi\left(\overline{e^{tA}}x+y\right)\,\caln(0,Q_{t})(dy)$ defined for $x\in \overline{H}$. It follows that the functions $g$ and $\overline{P}_{t}[\phi]$ are both continuous extensions of the function $P_{t}[\phi]$. Since the inclusion $H\hookrightarrow \overline{H}$ is dense and continuous we can conclude that $g(x)=\overline{P}_{t}[\phi](x)$ for all $x\in \overline{H}$.
\hfill\qed

\begin{remark}
From now on, motivated by the previous results, we will use the notation $\{P_{t}\}_{t\geq0}$ for both the Ornstein-Uhlenbeck semigroup and the extended one.
\end{remark}

Now, we recall some regularizing properties of the Ornstein-Uhlenbeck semigroup, see e.g. {\cite[Section 3]{LunardiPallara17}} and {\cite[Section 2, Theorem 2.1]{DaPratoZabczyk91}}.

\begin{hypothesis}
\label{ipo-smoothing}
\begin{itemize}
\item[]
\item [(i)]
For every $t> 0$, $h \in H$ we have
$e^{tA}h \in \operatorname{Im}
 Q_{t}^{1/2}$.
Consequently, by the closed graph Theorem, the operator
$\Lambda(t):H\to H$,
$\Lambda(t)h:=Q_{t}^{-1/2}
 e^{tA}h$,
for all $h \in H$,
is well defined for all $t>0$.
\vspace{-0.2truecm}
\item [(ii)]
There exists $\kappa_0>0$ and $\gamma \in (0,1)$ such that
$\|\Lambda(t)\|_{\call(H)} \le \kappa_0 \left(t^{-\gamma}\vee 1\right), \qquad
\forall t > 0$.
\end{itemize}
\end{hypothesis}

\begin{theorem}
Let $\phi\in B_{b}(H)$. Let $\{P_{t}\}_{t\geq0}$ be the Ornstein-Uhlenbeck semigroup defined on (\ref{ou-sgr}). Assume that \ref{ipo-smoothing} holds. Then $P_{t}[\phi]$ belongs to $C^{1}_{b}(H)$ and the following formula for the gradient holds
\begin{equation*}
\<\nabla P_{t}[\phi](x),h \>_{H}= \int_{H} \phi\left(e^{tA}x+y\right)\<\Lambda(t)h,Q_{t}^{-1/2}y\>_{H} \,\caln(0,Q_{t})(dy),
\end{equation*}\newline
Moreover, the following estimate is straightforward
\begin{equation*}
\left|\<\nabla P_{t}[\phi](x),h\>_{H}\right|\leq |h|_{H}\| \phi\|_{\infty} \|\Lambda(t)\|_{\call(H)}.
\end{equation*}
\end{theorem}

This theorem states the classical "strong Feller" property. It states that the semigroup smooths out measurable, bounded functions into continuously differentiable ones. The key condition $\operatorname{Im}(e^{tA})\subseteq \operatorname{Im}(Q_{t}^{1/2})$ is a form of null controllability condition for the linear system associated with the drift $A$ and diffusion $G$ {\cite[Section 9.4.6]{DaPratoZabczyk14}}, \cite[Part IV, Chapter 2]{Z92}.

\begin{remark}
The condition on the images $\operatorname{Im}(e^{tA})\subseteq \operatorname{Im}(Q_{t}^{1/2})$ in general is too strong and it is not satisfied in many applications, such as boundary control problems or stochastic equations with delay in the control. For this reason, we follow the ideas of \cite[Section 4.1]{FGFM-I}, \cite[Proposition 5.9]{FGFM-III}, \cite[Proposition 2.21]{GozziMasiero23} and we introduce a selection operator $P\in\call(H)$. The approach pioneered by Gozzi and Masiero is to seek a "partial" smoothing property that holds only in specific directions or for specific classes of functions, which is nonetheless sufficient to solve the control problem.
\end{remark} Hence, we define a special class of functions, depending on the the selection operator $P$.

\begin{definition}\label{defspecfun}
Let $P\in \call(H)$. We call $B_{b}^{P}(H)$ the subset of $B_{b}(H)$ of functions $\phi$ for which there exists a borel measurable and bounded function $\bar \phi: Im(P)\to\R$ such that
\begin{equation*}
\phi(x)=\bar \phi(P x), \quad x\in H.
\end{equation*}
\end{definition}

\begin{hypothesis}
\label{ipo-partial-smoothing}
\begin{itemize}
\item[]
\item [(i)]
For every $t> 0$, $k \in K$ we have
$Pe^{tA}Bk \in \operatorname{Im}
\left(P Q_{t} P^*\right)^{1/2}$.
Consequently, by the closed graph Theorem, the operator
$\Lambda^{P,B}(t):K\to H$,
$\Lambda^{P,B}(t)k:=\left(P Q_{t} P ^*\right)^{-1/2}
P e^{tA}Bk$,
for all $k \in K$,
is well defined for all $t>0$.
\vspace{-0.2truecm}
\item [(ii)]
There exists $\kappa_0>0$ and $\gamma \in (0,1)$ such that
$\|\Lambda^{P,B}(t)\|_{\call(K,H)} \le \kappa_0 \left(t^{-\gamma}\vee 1\right), \qquad
\forall t > 0$.
\end{itemize}
\end{hypothesis}
\begin{remark}
Hypothesis $(i)$ is the key controllability-like assumption for partial smoothing. It is a weaker version of the condition in the previous theorem, requiring the inclusion of images not in the whole space, but only after applying the selection operator $P$, while Hypothesis $(ii)$ controls the blow-up of the operator $\Lambda^{P,B}(t)$ as $t\to 0$. The condition $\gamma \in (0,1)$ ensures that the singularity is integrable, a crucial property for using fixed-point arguments to solve the HJB equation \cite[Hypothesis 5.7]{FGFM-III}, \cite[Hypothesis 2.19]{GozziMasiero23}.
\end{remark}
\begin{remark}\label{remark-ipo-partial smoothing-concreta}
By using a standard duality argument {\cite[Appendix B.2, Proposition B.1]{DaPratoZabczyk14}} we observe that \ref{ipo-partial-smoothing} $(i)$ is equivalent to the existence of a function $c(t)$ such that
\begin{equation}\label{eq:NCdualow}
\left|(Pe^{tA}B)^{*}z\right|^{2}_{K} \leq c(t) \< \left(Q_{t}P \right)^{*}z, P^{*}z\>_{H}, \quad t\in(0,\infty), \quad z\in H.
\end{equation}
Moreover, for any time $t>0$ the infimum over the constants $c(t)$ for which the inequality holds is exactly equal to the operator norm $\|\Lambda^{P,B}(t)\|_{\call(K,H)}$.
This dual formulation is often easier to verify in concrete applications. It relates the observability of the system in the direction of the control to the covariance of the noise, projected onto the selected subspace.
\end{remark}

The following is a generalization of the previous theorem, which recquires a weaker assumption on the images of the operators.
\begin{theorem}\label{teorema-reg-ou}
Let $\phi\in B_{b}^{P}(H)$. Let $\{P_{t}\}_{t\geq0}$ be the Ornstein-Uhlenbeck semigroup defined on (\ref{ou-sgr-ext}). Assume that $\operatorname{Im}(P e^{tA}B)\subseteq \operatorname{Im}((P Q_{t}P^{*})^{1/2})$ for some $t>0$. Then $P_{t}[\phi]$ belongs to $C^{1,B}_{b}(H)$ and the following formula for the $B$-gradient holds
\begin{equation*}
\<\nabla^{B}P_{t}[\phi](x),k\>_{K}= \int_{H} \phi\left(e^{tA}x+y\right)\<\Lambda^{P,B}(t)k,(P Q_{t} P^{*})^{-1/2}y\>_{H} \,\caln(0,Q_{t})(dy).
\end{equation*}\newline
Moreover, the following estimate is straightforward
\begin{equation*}
\left|\<\nabla^{B} P_{t}[\phi](x),k\>_{K}\right|\leq |k|_{K}\| \phi\|_{\infty} \|\Lambda^{P,B}(t)\|_{\call(K,H )}.
\end{equation*}
\end{theorem}

\subsection{"Lifted" partial smoothing}
In the previous section, we established a partial smoothing result for the Ornstein-Uhlenbeck semigroup for class of functions depending on a selection operator $P$. This result is general, but it is not clear how it can be applied to specific problems, such as problems with state dependent costs or with infinite time horizon. To overcome this limitation, we introduce here a more powerful and sophisticated technique, known as the "lifting map," first developed in \cite[Section 2]{GozziMasiero23}.
\\
The core idea of lifting is to move from a state-based perspective to a trajectory-based one. Instead of considering functions of the state $x \in \overline{H}$ at a single point in time, we will consider functions of the entire path $\{ \overline{Pe^{tA}}x \}_{t>0}$ generated by the state.

In this section, we will define this lifting map, denoted by $\Upsilon^P_\infty$, which maps a state $x$ to its corresponding path. We will then introduce a new class of functions, $\cals^P_\infty(\overline{H})$, which depend on the state only through its lifted trajectory. The main goal is to prove an extended partial smoothing result (Proposition \ref{prop-lift-partial-smoothing}) for this new class of functions, which will be the key tool to solve the HJB equation with state-dependent costs. Then, in what follows, we recall the results of \cite{GozziMasiero23} for the commodity of the reader.

Let $P:H\rightarrow H$ be a selection operator on $H$ as in the previous section. We make the following assumption throughout the discussion.

\begin{hypothesis}\label{ipo-traj}
For any $t>0$ the map $Pe^{tA}:H\rightarrow H$ can be extended into a linear and continuous map $\overline{Pe^{tA}}: \overline{H}\rightarrow H$.
\begin{remark}
This is a crucial regularity assumption. It requires that the composition of the evolution operator $e^{tA}$ with the selection operator $P$ is "smoothing" enough to map elements from the large, potentially rough space $\overline{H}$ into the well-behaved Hilbert space $H$. This property is satisfied in many relevant examples, including heat equations with boundary control, as shown in \cite[Section 4.1]{GozziMasiero23}.
\end{remark}
\end{hypothesis}

\begin{lemma}\label{lemma-ext}
Let Hypotheses \ref{ipo-data} and \ref{ipo-traj}
hold true.
For every $x\in \overline{H}$ and
$0<s\leq t$ we have
\vspace{-0.2truecm}
\begin{equation}\label{eq:Petasemigroup}
\overline{Pe^{tA}} x=
\overline{Pe^{sA}} \cdot \overline{e^{(t-s)A}} x.
\vspace{-0.2truecm}
\end{equation}
Moreover the map
$(0,+\infty)\to H$,
$t \mapsto \overline{Pe^{tA}} x$,
is continuous.
\end{lemma}

\begin{definition}\label{def-path-space}
We define the set of paths\vspace{-0.2truecm}
\[
\calc^P_A((0,\infty); H):= \left\{
f\in
C ((0,\infty); H) \hbox{ such that } \exists\, x\in \overline{H}:\;
f(t)=\overline{Pe^{tA}} x,\; \forall t \in (0,\infty)\right\}.\vspace{-0.2truecm}
\]
For every
$x\in \overline{H}$
we call
$y^P_{x}(\cdot)$
the path given by
$y^P_{x}(t)=\overline{Pe^{tA}}x$ for all $t \in (0,\infty)$.
This defines the space of all possible trajectories that can be generated by initial conditions in $\overline{H}$. This space will become the domain for the "lifted" functions we will work with.
\end{definition}

\begin{lemma}
\label{lemma-conv-path}
Let Hypotheses \ref{ipo-data} and \ref{ipo-traj}
hold true. Define the map
$\Upsilon^P_\infty:\overline{H} \to \calc^P_A((0,\infty); H)$
as $\Upsilon^P_\infty (x)=y^P_{x}$,
for all $x \in \overline{H}$.
$\Upsilon^P_\infty$ is surjective but not necessarily injective.
Moreover, we have, endowing $C((0,\infty); H)$ with the topology of uniform convergence on compact subsets, that
$x_n\to x$ in $\overline{H}$ $\Longrightarrow$ $y^P_{x_n}\to y^P_{x}$
in $C ((0,\infty); H)$.
The {converse} is not true in general.
Hence $\Upsilon^P_\infty$ is continuous if we endow $\calc^P_A((0,\infty); H)$ with the topology inherited by
$C ((0, \infty); H)$.
\begin{remark}
The map $\Upsilon^P_\infty$ is the "lifting map". It takes a point $x$ in the state space and maps it to a full trajectory. The continuity of this map ensures that small changes in the initial state lead to small changes in the resulting trajectory (in the appropriate topology). The lack of injectivity and the failure of the converse implication highlight that different initial states can lead to the same trajectory, and that convergence of trajectories does not imply convergence of the initial states. This is a key feature of infinite-dimensional systems \cite[Lemma 2.10]{GozziMasiero23}.
\end{remark}
\end{lemma}

\begin{remark}
Without explicit notice we will take,
on $\mathcal{C}^P_A((0,\infty); H)$, the topology inherited by
$C ((0,\infty); H)$.
\end{remark}

\begin{definition}\label{def-L2-path}
We define the space of functions that are square integrable on $[0,+\infty)$ when multiplied by a suitable weight $e^{-\rho t}, \,\rho>0$:\vspace{-0.2truecm}
\[
L^2_\rho(0,\infty;H):=\left\{ f:[0,+\infty)\to H: \quad \|f\|_{L_\rho^2}:=
\left( \int_0^{+\infty}e^{-2\rho t} \left| f(t)\right|^2_H\,dt \right)^{1/2}<+\infty \right\}. \vspace{-0.2truecm}
\]
This is a Hilbert space with the inner product
$\<f,g\>_{L_\rho^2}:= \int_0^{+\infty}
e^{-2\rho t} \<f(t),g(t)\>_H\,dt$.
For simplicity, if no confusion is possible we will write $L^2_\rho$ instead of $L^2_\rho(0,\infty;H)$ and $L^2_T$ for $L^2(0,T;H)$.
\end{definition}

\begin{hypothesis}
\label{ipo-growth-paths}
There exist $\omega\in\R$, $C>0$ and $\eta \in [0,1/2)$ such that
$\left|\overline{Pe^{tA}} x\right|_H\leq e^{\omega t}t^{-\eta}\left|x\right|_{\overline{H}}$, $\forall x \in \overline{H}$,
and, as a consequence, the map
$t \to \overline{Pe^{tA}} x$,
belongs to $L^2_{\rho}(0,\infty;H)$.
This assumption provides a quantitative bound on the growth and singularity of the paths. The term $t^{-\eta}$ allows for a possible blow-up at $t=0$, which is typical for heat-like semigroups, while $e^{\omega t}$ controls the growth at infinity. The condition $\eta < 1/2$ is crucial to ensure that the paths are square-integrable near the origin.
\end{hypothesis}
\begin{lemma}
\label{lm:UpsilonL2}
Let Hypotheses \ref{ipo-data}, \ref{ipo-traj}
and \ref{ipo-growth-paths} hold true. Let $A$ be of type $\omega${\footnote{We say that $A$ is of type $\omega$ if there exists some $M\geq 1$ such that $\left\| e^{At}\right\| \leq Me^{\omega t}$ for all $t\ge 0$.}} and consider $\rho>\omega$. then for all $x\in \overline{H}$, we have
$ y^P_{x}\in L_\rho^2(0,+\infty; H)$.
Hence,
$\calc^P_A((0,\infty); H)$
can be seen as a linear subspace of $L_\rho^2(0,\infty; H)$ with continuous embedding given by
$\Upsilon^P_{\infty}$.
Moreover, if $x_n\to x$ in $\overline{H}$ this implies that $y^P_{x_n}\to y^P_{x}$ in $L_\rho^2 (0,\infty; H)$ but the {converse} is not true in general.
\end{lemma}

\begin{remark}\label{rem:space}
The extension to the infinite horizon requires a careful choice of the trajectory space. In the finite horizon case (as in \cite{GozziMasiero23}), trajectories are naturally defined on a compact interval $[0,T]$, and spaces like $L^2(0,T;H)$ are sufficient. In our stationary case, we must consider trajectories over $(0,\infty)$. To ensure that the trajectories $y_x^P(\cdot)$ belong to a well-defined Hilbert space (as required by Hypothesis \ref{ipo-growth-paths} and Lemma \ref{lm:UpsilonL2}) and to manage their potential growth at infinity, we are forced to introduce the weighted space $L^2_\rho(0,\infty;H)$ (Definition \ref{def-L2-path}). The weight $e^{-\rho t}$ is essential to ensure integrability and the well-posedness of the lifting operator $\Upsilon^P_\infty$, making the technical setup more complex compared to the unweighted case on a finite interval.
\end{remark}

\begin{definition}
\label{df:SP}
Let $L$ be another separable Banach space. We introduce the following set of functions
\begin{equation}\label{def-funz-traj}
\cals^P_\infty(\overline{H};L)
:=\left\{\phi: \overline{H}\to L \hbox{ s.t. }
\exists \,\hat \phi:\calc^P_A((0,\infty); H) \to L :
\hbox{ $\hat \phi$ bounded,
Borel meas. and }
\phi(x)=\hat \phi\left(y^P_{x}\right), \;\forall x \in \overline{H}.
\right\}. \end{equation}
When $L$ coincides with $\R$ we will write directly $S^{P}_{\infty}(\overline{H})$.
\end{definition}
This is the central definition of this section. It introduces the class of "lifted" functions $\cals^P_\infty(\overline{H})$. A function $\phi$ belongs to this class if its value at a point $x$ depends only on the entire trajectory $y^P_x$ originating from $x$.

\begin{proposition}\label{prop-adjoint}
Assume Hypotheses \ref{ipo-data}, \ref{ipo-traj} and \ref{ipo-growth-paths}.
We have the following.
\begin{enumerate}
\item [(i)]
We have $\cals^P_{\infty}(\overline{H})
\subseteq B_b(\overline{H})$. Moreover the restriction to $H$ of a function in $\cals^P_{\infty}(\overline{H})$ belongs to
$B_b(H)$.
\item [(ii)]
Let $\phi\in \cals^P_\infty(\overline{H})$ and $\hat\phi$
be the function given in the definition of $\cals^P_\infty(\overline{H})$: we have
$\phi=\hat \phi \circ \Upsilon^P_\infty$.
\item [(iii)]
Assume that $P$ can be extended to a continuous linear operator $\overline{P}:\overline{H}\to H$ such that
$\operatorname{Im} \overline{P}=\operatorname{Im} P$. Then
we have
$B_b^P(H)\subseteq\cals^P_\infty(H)$. Moreover, if also $P$ commutes with $A$,
then $B_b^P(H)=\cals^P_\infty(H)$.
\item [(iv)]
The adjoint operator
$(\Upsilon^P_\infty)^*:L^2_\rho(0,\infty;H)
\to \overline{H}'\subseteq H$
is given by
$(\Upsilon^P_\infty)^*z(\cdot)=\int_{0}^{+\infty}e^{-\rho s}e^{sA^*}P^*z(s)\,ds$.
\end{enumerate}

\end{proposition}

This proposition establishes the fundamental properties of the class $\cals^P_\infty(\overline{H})$. Point (iii) clarifies the relationship between this new class and the more classical function classes depending only on $Px$, showing that the lifting approach is a true generalization. Point (iv) provides an explicit formula for the adjoint of the lifting operator, which will be essential for verifying the smoothing hypothesis (Hypothesis \ref{ipo-smoothing-lifting}) later on.

\begin{hypothesis}
\label{ipo-smoothing-lifting}
\begin{itemize}
%\item[]
\item [(i)]
For every $t> 0$, $k \in K$ we have
$\Upsilon^P_\infty\overline{e^{tA}}Bk \in \operatorname{Im}
\left(\Upsilon^P_\infty Q_{t} (\Upsilon^P_\infty)^*\right)^{1/2}$.
Consequently, by the closed graph Theorem, the operator
$\widehat\Lambda^{P,B}(t):K\to L^2_\rho(0,\infty;H)$,
$\widehat\Lambda^{P,B}(t)k:=\left(\Upsilon^P_\infty Q_{t} (\Upsilon^P_\infty)^*\right)^{-1/2}
\Upsilon^P_\infty e^{tA}Bk$,
for all $k \in K$,
is well defined for all $t>0$.
\vspace{-0.2truecm}
\item [(ii)]
There exists $\kappa_0>0$ and $\gamma \in (0,1)$ such that
$\|\widehat\Lambda^{P,B}(t)\|_{\call(K,L^2_\rho)} \le \kappa_0 \left(t^{-\gamma}\vee 1\right), \qquad
\forall t > 0$.
\end{itemize}
This is the "lifted" version of the partial smoothing hypothesis (compare with Hypothesis \ref{ipo-partial-smoothing}).
\end{hypothesis}
This is the technical heart of the lifting method. It states that the trajectory generated by the control action, $\Upsilon^P_\infty\overline{e^{tA}}Bk$, is "less singular" than the noise projected onto the trajectory space, which is measured by its covariance operator $\Upsilon^P_\infty Q_{t} (\Upsilon^P_\infty)^*$. Verifying this hypothesis is the main challenge when applying the theory to a concrete example.

\begin{remark}\label{remark commutatività}
Following \cite[Remark 2.20]{GozziMasiero23} it can be verified that if the operator $P$ and the semigroup $e^{tA}$ commute for each $t\geq0$ then hypothesis (\ref{ipo-smoothing-lifting}) is equivalent to (\ref{ipo-partial-smoothing}), which is usually easier to verify in applications as we will see in our examples.
\end{remark}

\begin{remark}\label{remark-ipo-concreta}
Again, by using {\cite[Appendix B.2, Proposition B.1]{DaPratoZabczyk14}} and the fact that inclusion $H\hookrightarrow \overline{H}$ is dense and continuous we can conclude that Hypothesis \ref{ipo-smoothing-lifting} $(i)$ is equivalent to the existence of a function $c(t)$ such that
\begin{equation}\label{eq:NCdualhigh}
\left|(\Upsilon^P_\infty e^{tA}B)^{*}z\right|^{2}_{K} \leq c(t) \< (Q_{t}\Upsilon^P_\infty)^{*}z, (\Upsilon^P_\infty)^{*}z\>_{H}, \quad t\in(0,\infty), \quad z\in H.
\end{equation}
Moreover, for any time $t>0$ the infimum over the constants $c(t)$ for which the inequality holds is exactly equal to the operator norm $\|\widehat\Lambda^{P,B}(t)\|_{\call(K,L^2_\rho)}$.
\end{remark}

\begin{proposition}\label{prop-lift-partial-smoothing}
Let Hypotheses \ref{ipo-data}, \ref{ipo-traj}, \ref{ipo-growth-paths}, and
\ref{ipo-smoothing-lifting}-(i)
hold true.
Then $P_t,\,t>0$ maps functions $\phi\in \cals^P_\infty(\overline{H})$
into functions which are $B$-Fr\'echet differentiable in $\overline{H}$, and the $B$-derivative is given,
for all $t>0$, $x\in \overline{H}$, by\vspace{-0.2truecm}
\begin{align}\label{eq:formulader-gen-Pnew}
&\nabla^B(P_{t}[\phi])(x)k =
\int_{L^2_\rho}\hat\phi\left(z_1+ \Upsilon^P_\infty x \right)
\<\widehat\Lambda^{P,B}(t) k, (\Upsilon^P_\infty Q_t(\Upsilon^P_\infty)^*)^{-1/2} z_1\>_{L^2_\rho}
\,\caln(0,(\Upsilon^P_\infty Q_t(\Upsilon^P_\infty)^*))(dz_1)
\notag
\\[2mm]
&=
\E\left[\hat\phi\left(\Upsilon^P_\infty
X(t; x)\right)
\<\widehat\Lambda^{P,B}(t) k,
(\Upsilon^P_\infty Q_t(\Upsilon^P_\infty)^*)^{-1/2}\Upsilon^P_\infty W_A(t)\>_{L^2_\rho}
\right].
\end{align}
Moreover, for any $\phi\in \cals^P_\infty(\overline{H})$, $t>0$, $x\in \overline{H}$, $k\in K$,\vspace{-0.2truecm}
\begin{equation}\label{norm-Cdernew}
\left| \<\nabla^B P_t[\phi](x), k\>\right| \leq \|\widehat\Lambda^{P,B}(t)\|_{\call(K,L^2_\rho)}
\Vert\phi\Vert_\infty \cdot\left| k\right|.
\end{equation}
\end{proposition}
This proposition is the culmination of the lifting strategy. It establishes that the O-U semigroup $P_t$ regularizes functions from the lifted class $\cals^P_\infty(\overline{H})$ to make them differentiable in the direction of the control operator $B$. The formula for the derivative and the associated estimate are the essential ingredients for solving the stationary HJB equation with state-dependent running costs via a fixed-point argument.

\section{Solution of the Hamilton-Jacobi-Bellman Equation}\label{sec:hjb_equation}
In this section, we arrive at the core of the optimal control problem: solving the associated Hamilton-Jacobi-Bellman (HJB) equation. The HJB equation provides the link between the value function $V(x)$ and the dynamics of the system. Thanks to the "lifting" technique and the partial smoothing results established in the previous section, we are now equipped to handle stationary HJB equations with state-dependent costs, which was our main objective.
\\
We will first formally write down the stationary HJB equation for our infinite-horizon problem. Then, we will define what we mean by a "mild solution", which is an integral formulation of the HJB equation that is well-suited for our setting. The main result of this section, Theorem \ref{teorema-esistenza-fixed}, will be to prove the existence and uniqueness of such a mild solution by means of a fixed-point argument on a carefully constructed operator. This result relies heavily on the regularizing properties of the O-U semigroup for lifted functions. Finally, we will investigate the higher-order regularity of this solution.

\subsection{The HJB Equation and Mild Solutions}

We define the current value Hamiltonian for $p,u\in K$
\[
H_{CV}(p\,;u):=\<p,u\>_{K}+\ell_1(u)
\]
and the minimum value Hamiltonian
\begin{equation}\label{psi1-gen}
H_{min}(p)=\inf_{u\in U}H_{CV}(p\,;u).
\end{equation}
The HJB equation associated to the stochastic optimal control problem is formally
\begin{equation}\label{HJB-ellittica}
\lambda v(x)=\cala [v(\cdot)](x) +\ell_0(x)+
H_{min}\left(\nabla^B v(x)\right), \quad x\in H.
\end{equation}
Here the differential operator $\cala$
is the infinitesimal generator of the O-U
semigroup $(P_{t})_{t\geq 0}$, formally defined by
\begin{equation}\label{eq:ell-gen}
\cala[f](x)=\frac{1}{2} \operatorname{Tr} \left(Q\nabla^2f(x)\right)
+ \< x,A^*\nabla f(x)\>_{H}.
\end{equation}
\begin{remark}
This is the stationary HJB equation associated with the infinite horizon problem. The term $\lambda v(x)$ comes from the discounting in the cost functional. The operator $\mathcal{A}$ represents the expected change of the value function along the uncontrolled dynamics, while the term $H_{min}(\nabla^B v(x))$ is the nonlinear part that encodes the optimization over the controls. Note that the Hamiltonian depends on the generalized $B$-derivative $\nabla^Bv$, which is why the partial smoothing theory is so critical.
\end{remark}

\begin{definition}\label{def-mild-HJB}
We say that a function $v:\overline{H}\rightarrow\R$ is a mild solution of (\ref{HJB-ellittica}) if
\begin{itemize}
\item[(i)] $v\in C^{1,B}_{b}(\overline{H})$;
\item[(ii)] for each $x\in \overline{H}$ it satisfies the following integral equation
\begin{equation}\label{hjb-mild}
v(x)= \int_{0}^{\infty} e^{-\lambda t} P_{t}\left[\ell_{0}(\cdot)+H_{min}\left(\nabla^{B}v(\cdot)\right)\right](x)\,dt
\end{equation}
\end{itemize}
\end{definition}
This definition is based on the variation of constants formula. Instead of satisfying the PDE in a classical sense, the solution satisfies an equivalent integral equation. This approach avoids issues with the domain of the operator $\mathcal{A}$ and is standard for solving HJB equations in infinite dimensions {\cite[Section 3, Definition 3.1]{FedericoGozzi17}}{\cite[Chapter 13, Section 13.2.1]{DaPratoZabczyk02}}.

\begin{remark}
Since in our framework the semigroup $P_{t}$ can be naturally extended to a semigroup acting on functions of $B_{b}(\overline{H})$ and the running cost $\ell_{0}$ is defined on $\overline{H}$, equation (\ref{hjb-mild}) can be studied in the extended space $\overline{H}$. Hence, form now on, our goal will be to study equation (\ref{HJB-ellittica}) in mild form on the extended state space $\overline{H}$.
\end{remark}

\begin{definition}\label{def-operatore-integrale}
We introduce the following nonlinear operator $F=(F_{1},F_{2})$, acting on $S^{P}_{\infty}(\overline{H})\times S^{P}_{\infty}(\overline{H};K)$ by
\begin{equation*}
F_{1}[v,w](x)= \int_{0}^{\infty} e^{-\lambda t} P_{t}\left[\ell_{0}(\cdot)+H_{min}\left(w(\cdot)\right)\right](x)\,dt
\end{equation*}
\begin{equation*}
F_{2}[v,w](x)= \int_{0}^{\infty} e^{-\lambda t} \nabla^{B}P_{t}\left[\ell_{0}(\cdot)+H_{min}\left(w(\cdot)\right)\right](x)\,dt
\end{equation*}
\end{definition}
Solving the mild HJB equation \eqref{hjb-mild} is equivalent to finding a fixed point for the operator $F_1$. The operator is defined on the space of lifted functions $\cals^P_\infty(\overline{H})$ because the presence of the state-dependent cost $\ell_0$ means we must leverage the smoothing properties on this specific class of functions, as established by Proposition \ref{prop-lift-partial-smoothing}.

\begin{remark}\label{remark-ben-def}
The operator $F$ is well defined due to the regularizing property of the O-U semigroup discussed in (\ref{prop-lift-partial-smoothing}). Infact for any function $w \in S^{P}_{\infty}(\overline{H};K)$ the function $\ell_{0}+H_{min}(w)$ belongs to $S^{P}_{\infty}(\overline{H})$ and so it is $B$-regularized by the O-U semigroup. Moreover, it is immediate to check that $F$ takes value in $S^{P}_{\infty}(\overline{H})\times S^{P}_{\infty}(\overline{H};K)$.
\end{remark}
This is the main result of this section. It establishes existence and uniqueness of the HJB solution by applying the Banach fixed-point theorem. The proof will show that for a sufficiently large discount factor $\lambda$, the operator $F$ is a contraction on the space of lifted functions. The magnitude of $\lambda$ is needed to absorb the Lipschitz constants and operator norms that appear in the estimates.
\begin{theorem}\label{teorema-esistenza-fixed}
Assume that hypothesis (\ref{ipo-data}), (\ref{ipo-traj}), (\ref{ipo-growth-paths}) and (\ref{ipo-smoothing-lifting}) holds. Then, there exist a $\lambda_{0}>0$ such that for each $\lambda\geq\lambda_{0}$ the nonlinear operator $F$ introduced in (\ref{def-operatore-integrale}) admits a unique fixed point in the space $S^{P}_{\infty}(\overline{H})\times S^{P}_{\infty}(\overline{H};K)$.
\end{theorem}

\dim
As explained in Remark (\ref{remark-ben-def}) the operator $F$ maps $S^{P}_{\infty}(\overline{H})\times S^{P}_{\infty}(\overline{H};K)$ into itself. We show that $F$ is actually a contraction in this space for sufficiently big discount factor $\lambda>0$. Let $x\in \overline{H}$, $v,\tilde v\in S^{P}_{\infty}(\overline{H}) $ and $w, \tilde w \in S^{P}_{\infty}(\overline{H};K) $. For the first component of $F$ we have
\begin{align*}
& \left|F_{1}[v,w](x)-F_{1}[\tilde v, \tilde w](x)\right|= \left|\int_{0}^{\infty}e^{-\lambda t}P_{t}\left[H_{min}(w)-H_{min}(\tilde w)\right](x)\,dt \right|\\
& \leq \left|\int_{0}^{\infty}e^{-\lambda t}\int_{\overline{H}}\left[H_{min}\left(w(y)\right)-H_{min}\left(\tilde w(y)\right)\right]\,\caln\left(\overline{e^{tA}}x,Q_{t}\right)(dy)\,dt \right| \\
& \leq \int_{0}^{\infty}e^{-\lambda t}\int_{\overline{H}}\left|H_{min}\left(w(y)\right)-H_{min}\left(\tilde w(y)\right)\right|\,\caln\left(\overline{e^{tA}}x,Q_{t}\right)(dy)\,dt.
\end{align*}
By using the lipschitzianity of $H_{min}$ we can estimate the last term with
\begin{align*}
C\int_{0}^{\infty}e^{-\lambda t}\int_{H}\left|w(y)-\tilde w(y)\right|\,\caln\left(\overline{e^{tA}}x,Q_{t}\right)(dy)\,dt\leq C\|w-\tilde w\|_{\infty}\int_{0}^{\infty} e^{-\lambda t} dt.
\end{align*}
For the second component of operator $F$, similarly, we have
\begin{align*}
& \left|F_{2}[v,w](x)-F_{2}[\tilde v, \tilde w](x)\right|= \left|\int_{0}^{\infty}e^{-\lambda t}\nabla ^{B}P_{t}\left[H_{min}(w)-H_{min}(\tilde w)\right](x)\,dt \right|\\
& \leq \int_{0}^{\infty}e^{-\lambda t}\left|\nabla ^{B}P_{t}\left[H_{min}(w)-H_{min}(\tilde w)\right](x)\right|\,dt.
\end{align*}
Thanks to Proposition (\ref{prop-lift-partial-smoothing}) we can estimate
\begin{equation*}
\int_{0}^{\infty}e^{-\lambda t}\left|\nabla ^{B}P_{t}\left[H_{min}(w)-H_{min}(\tilde w)\right](x)\right|\,dt \leq \int_{0}^{\infty}e^{-\lambda t} \left\| H_{min}(w)-H_{min}(\tilde w)\right\|_{\infty} \|\widehat\Lambda^{P,B}(t)\|_{\call(K,L^2_\rho)}\,dt.
\end{equation*}
By using the lipschitzianity of $H_{min}$ and Hypothesis (\ref{ipo-smoothing-lifting}) we have
\begin{equation*}
\int_{0}^{\infty}e^{-\lambda t} \left\| H_{min}(w)-H_{min}(\tilde w)\right\|_{\infty} \|\widehat\Lambda^{P,B}(t)\|_{\call(K,L^2_\rho)}\,dt \leq C\|w-\tilde w\|_{\infty}\int_{0}^{\infty}e^{-\lambda t} \left(1 \vee t^{-\gamma}\right)\,dt
\end{equation*}
for a certain constant $C>0$ and $\gamma\in(0,1)$. Putting all together we obtain
\begin{align*}
\left| F[v,w](x)-F[\tilde v,\tilde w](x)\right|\leq C \|w-\tilde w\|_{\infty}\int_{0}^{\infty} e^{-\lambda t} \left(1 \vee t^{-\gamma}\right)\,dt.
\end{align*}
It is easy to check that the function $\lambda\to\int_{0}^{\infty} e^{-\lambda t} \left(1 \vee t^{-\gamma}\right)\,dt$ is monotone decreasing and goes to $0$ as $\lambda\to\infty$. Hence, there exists a $\lambda_{0}>0$ such that for each $\lambda\geq\lambda_{0}$ the term $C\int_{0}^{\infty} e^{-\lambda t} (1 \vee t^{-\gamma})\,dt$ is strictly less then one, showing that $F$ is a contraction. By the contraction mapping theorem $F$ admits a unique fixed point.
\qed

\begin{remark}\label{rem:fixedpoint}When solving the evolutionary (finite horizon) HJB equation via a fixed-point argument, the contraction is often obtained by integrating over a sufficiently small time interval $[t,T]$. This is not possible in the stationary case, where the integral operator (Definition \ref{def-operatore-integrale}) is defined over $[0,\infty)$. As shown in the proof of Theorem \ref{teorema-esistenza-fixed}, obtaining a contraction requires the discount factor $\lambda$ to be large.
We are currently working on an extension of our results to remove such restriction.\end{remark}

\begin{lemma}\label{lemma-mild-fixed}
Let $(v,w)\in S^{P}_{\infty}(\overline{H})\times S^{P}_{\infty}(\overline{H};K)$ be the unique fixed point of the operator $F$ introduced in (\ref{def-operatore-integrale}). Then, one has
\begin{itemize}
\item[(i)] $v\in C^{1,B}(\overline{H})$;
\item[(ii)] $\nabla^{B}v=w$;
\item[(iii)] $v$ is the unique mild solution of equation (\ref{HJB-ellittica}) in the class $S^{P}_{\infty}(\overline{H})$.
\end{itemize}
\end{lemma}

\dim
By construction the function $v$ satisfies
\begin{equation}\label{1-2-comp}
v(x)=\int_{0}^{\infty} e^{-\lambda t} P_{t}\left[\ell_{0}(\cdot)+H_{min}\left(w(\cdot)\right)\right](x)\,dt, \quad x\in \overline{H}.
\end{equation}
By using proposition (\ref{prop-lift-partial-smoothing}) we get that $v\in C^{1,B}_{b}(\overline{H})$ and we can differentiate under the integral sign obtaining
\begin{equation*}
\nabla^{B}v(x)=\int_{0}^{\infty} e^{-\lambda t}\nabla^{B} P_{t}\left[\ell_{0}(\cdot)+H_{min}\left(w(\cdot)\right)\right](x)\,dt,
\end{equation*}
with the last term coinciding with $w(x)$ with $x\in \overline{H}$, proving (i) and (ii). Now, we can substitute $w$ with $\nabla ^{B}v$ in equation (\ref{1-2-comp}) to get
\begin{equation*}
v(x)=\int_{0}^{\infty} e^{-\lambda t} P_{t}\left[\ell_{0}(\cdot)+H_{min}\left(\nabla^{B}v(\cdot)\right)\right](x)\,dt, \quad x\in \overline{H},
\end{equation*}
which conclude the proof of (iii).
\qed

\begin{remark}
The previous theorem clearly ensures the existence of mild solutions of (\ref{HJB-ellittica}) at least for large values of the discount factor $\lambda>0$. However, uniqueness is guaranteed only on the class $S^{P}_{\infty}(\overline{H})$.
\end{remark}

\subsection{Higher-Order Regularity}
The following are additional regularity assumptions on the problem data, which will allow us to prove higher regularity for the solution of the HJB equation.
\begin{hypothesis}\label{ipo-ham-reg}
\begin{itemize}
\item [(i)] $\ell_{0}:\overline{H}\to\R$ is of class $C^{1}_{b}(\overline{H})$,
\item [(ii)] $H_{min}:K\to\R$ is of class $C^{1}_{b}(K)$,
\item [(iii)] there exists a constant $L>0$ such that
\begin{equation*}
\left|\nabla\ell_{0}(x)-\nabla \ell_{0}(y)\right|+\left|\nabla H_{min}(p)- \nabla H_{min}(q)\right|\leq L\left(\left|x-y\right|+\left|p-q\right|\right)
\end{equation*}
uniformly for $x,y\in \overline{H}$ and $p,q\in K$.
\end{itemize}
\end{hypothesis}
This proposition shows that if we assume more regularity on the cost functions, we obtain more regularity on the value function. Proving that the solution is of class $C^{2,B}$ (i.e., that its $B$-derivative is itself differentiable) is a key step towards establishing the existence of optimal feedback controls via a verification theorem \cite{FGFM-II}
\begin{proposition}\label{prop-regolarità}
Assume that hypothesis (\ref{ipo-data}), (\ref{ipo-traj}), (\ref{ipo-growth-paths}) and (\ref{ipo-smoothing-lifting}) hold. Assume, morover, that hypothesis (\ref{ipo-ham-reg}) is satisfied. Let $\lambda_{0}>0$ as in Theorem (\ref{teorema-esistenza-fixed}). Then there exists $\lambda_{1}\geq \lambda_{0}$ such that for any $\lambda\geq\lambda_{1}$ the corresponding solution $v:\overline{H}\to\R$ of equation (\ref{HJB-ellittica}) is of class $C^{2,B}_{b}(\overline{H})$.
\end{proposition}
\dim
We introduce the Banach space
\begin{equation*}
\calh:=S^{1,P}_{\infty}(\overline{H})\times S^{1,P}_{\infty}(\overline{H};K)
\end{equation*}
equipped with the norm
\begin{equation*}
\|(u,w)\|_{\calh}:= \|(u,w)\|_{\infty}+ \|(\nabla u,\nabla w)\|_{\infty}.
\end{equation*}
Let $B_{1}$ be the closed unit ball in $\calh$. We want to prove that $F$ is a contraction on $B_{1}$ if $\lambda>0$ is sufficiently large. Let $(v,w)\in B_{1}$. We need to show that $F[v,w]$ still belongs to $\calh$ and in particular to $B_{1}$. For any $w\in S^P_{\infty}(\overline{H};K)$ we introduce the function $\phi^{w}:\overline{H}\to\R$ defined by
\begin{equation*}
\phi^{w}(x)=\ell_{0}(x)+H_{min}\left(w(x)\right), \quad x\in \overline{H}.
\end{equation*}
Since both $\ell_{0}$ and $H_{min}(w)$ belong to $S^{P}_{\infty}(\overline{H})$ there exists a unique measurable function $\hat \phi^{w}:\calc^P_A((0,\infty); H)\to\R$ such that $\phi^{w}(x)=\hat \phi^{w}\left(\Upsilon^P_\infty x\right)$ for all $x\in \overline{H}$.
Hence, the first component of $F$ is given by
\begin{equation}\label{eqf1}
F_{1}\left[v,w\right]\left(x\right)=\int_{0}^{\infty} e^{-\lambda t}P_{t}\left[\phi^{w}\left(\cdot\right)\right]\left(x\right)\,dt= \int_{0}^{\infty} e^{-\lambda t}\int_{H}\phi^{w}\left(\overline{e^{tA}}x+y\right)\,\caln(0,Q_{t})(dy)\,dt.
\end{equation}

By using hypothesis (\ref{ipo-ham-reg}) we can differentiate inside the integral the first component of $F$, obtaining the following
\begin{equation}\label{eqnf1}
\<\nabla F_{1}[v,w](x),h\>_{\overline{H}',\overline{H}}= \int_{0}^{\infty} e^{-\lambda t}\int_{H}\<\nabla\phi^{w}\left(\overline{e^{tA}}x+y\right), \overline{e^{tA}}h\>_{\overline{H}',\overline{H}}\,\caln(0,Q_{t})(dy)\,dt, \quad x,h\in \overline{H}.
\end{equation}
Now, we can write explicitely $\<F_{2}[v,w](x),k\>_{K}$ for each $k\in K$ using (\ref{prop-lift-partial-smoothing}) as
\begin{equation}\label{eqf2}
\int_{0}^{\infty} e^{-\lambda t} \int_{L^2_\rho}\hat\phi^{w}\left(z_1+ \Upsilon^P_\infty x \right)
\<\widehat\Lambda^{P,B}(t) k, (\Upsilon^P_\infty Q_t(\Upsilon^P_\infty)^*)^{-1/2} z_1\>_{L^2_\rho}
\,\caln(0,(\Upsilon^P_\infty Q_t(\Upsilon^P_\infty)^*))(dz_1)\,dt.
\end{equation}
Hence, by differentiating, for each $h\in \overline{H}$ and $k\in K$ we can compute $\<\nabla F_{2}[v,w](x)h,k\>_{K}$ obtaining
\begin{equation}\label{eqnf2}
\int_{0}^{\infty} e^{-\lambda t} \int_{L^2_\rho}\<\nabla \hat\phi^{w}\left(z_1+ \Upsilon^P_\infty x \right),\Upsilon^P_\infty h\>_{L^{2}_{\rho}}
\<\widehat\Lambda^{P,B}(t) k, (\Upsilon^P_\infty Q_t(\Upsilon^P_\infty)^*)^{-1/2} z_1\>_{L^2_\rho}
\,\caln(0,(\Upsilon^P_\infty Q_t(\Upsilon^P_\infty)^*))(dz_1)\,dt.
\end{equation}
By using equation (\ref{eqf1}) and the global lipschitzianity of $\ell_{0}$ and $H_{min}$ we can estimate
\begin{align}\label{estf1}
& \left|F_{1}[v,w](x)-F_{1}[\tilde v, \tilde w](x)\right|\leq \int_{0}^{\infty} e^{-\lambda t}\int_{H}\left|\phi^{w}\left(\overline{e^{tA}}x+y\right)-\phi^{\tilde w}\left(\overline{e^{tA}}x+y\right)\right|\,\caln(0,Q_{t})(dy)\,dt\leq\nonumber\\
& \leq \left\| \phi^{w}-\phi^{\tilde w} \right\|_{\infty} \int_{0}^{\infty} e^{-\lambda t}\,dt\leq \|w-\tilde w \|_{\infty} \int_{0}^{\infty} e^{-\lambda t}\,dt.
\end{align}
In the same way, by using equation (\ref{eqf2}) and hypothesis (\ref{ipo-smoothing-lifting}) for $|k|_{K}=1$ one can estimate
\begin{align}\label{estf2}
&\left|\<F_{2}[v,w](x)-F_{2}[\tilde v, \tilde w](x),k \>_{K}\right| \leq \left\| \phi^{w}-\phi^{\tilde w} \right\|_{\infty} \int_{0}^{\infty} e^{-\lambda t}\|\widehat\Lambda^{P,B}(t)\|_{\call(K,L^2_\rho)}\,dt \leq\nonumber\\
&\leq C \|w-\tilde w \|_{\infty} \int_{0}^{\infty} e^{-\lambda t}\left(1 \vee t^{-\gamma}\right)\,dt,
\end{align}
for some constant $C>0$ and $\gamma\in(0,1)$. Putting together these two estimates we deduce that
\begin{equation}
\|F_{1}[v,w]-F_{1}[\tilde v, \tilde w]\|_{\infty}+\|F_{2}[v,w]-F_{2}[\tilde v, \tilde w]\|_{\infty}\leq C \|w-\tilde w \|_{\infty} \int_{0}^{\infty} e^{-\lambda t}\left(1 \vee t^{-\gamma}\right)\,dt.
\end{equation}

By using equation (\ref{eqnf1}) we obtain for any unit vector $h\in \overline{H}$
\begin{align*}
& \left|\<\nabla F_{1}[v,w](x)-\nabla F_{1}[\tilde v, \tilde w](x),h\>_{\overline{H}',\overline{H}} \right| \leq \left\| \nabla\phi^{w}-\nabla\phi^{\tilde w} \right\|_{\infty} \int_{0}^{\infty} e^{-\lambda t}\|\overline{e^{tA}}\|_{\call(\overline{H})}\,dt \leq \\
& \leq \left\| \nabla\phi^{w}-\nabla\phi^{\tilde w} \right\|_{\infty} \int_{0}^{\infty} e^{-(\lambda-\omega)t}\,dt,
\end{align*}
being $\omega$ the type of the semigroup $\{\overline{e^{tA}}\}_{t\geq0}$.
Similarly, using equation (\ref{eqnf2}) we get
\begin{align}\label{estnf2}
&\left|\<\nabla F_{2}[v,w](x)-\nabla F_{2}[\tilde v, \tilde w](x),k \>_{K}\right| \leq C \left\| \nabla\phi^{w}-\nabla \phi^{\tilde w} \right\|_{\infty} \int_{0}^{\infty} e^{-\lambda t}\|\widehat\Lambda^{P,B}(t)\|_{\call(K,L^2_\rho)}\,dt \leq\nonumber\\
&\leq C \left\| \nabla\phi^{w}-\nabla \phi^{\tilde w}\right\|_{\infty} \int_{0}^{\infty} e^{-\lambda t}\left(1 \vee t^{-\gamma}\right)\,dt.
\end{align}
The last two estimates imply
\begin{equation}\label{estnf12}
\|\nabla F_{1}[v,w]-\nabla F_{1}[\tilde v,\tilde w]\|_{\infty}+\|\nabla F_{2}[v,w]-\nabla F_{2}[\tilde v,\tilde w]\|_{\infty}\leq C \left\| \nabla\phi^{w}-\nabla \phi^{\tilde w}\right\|_{\infty} \int_{0}^{\infty} e^{-(\lambda-\omega)t}\left(1 \vee t^{-\gamma}\right)\,dt.
\end{equation}
Moreover, it is easy to check that
\begin{equation}\label{estnfw}
\left\| \nabla\phi^{w}-\nabla \phi^{\tilde w}\right\|_{\infty} \leq C \left(\|\nabla w-\nabla \tilde w\|_\infty + \|w\|_\infty \|w-\tilde w\|_{\infty}\right)\leq C\|w-\tilde w\|_{\calh}
\end{equation}
and the last estimate follows from the fact that the functional is restricted to the unit ball $B_{1}$.
Combining (\ref{estnf12}) and (\ref{estnfw}) we get
\begin{equation}
\| F[v,w]- F[\tilde v,\tilde w]\|_{\calh}\leq C \| w-\tilde w\|_\calh \int_{0}^{\infty} e^{-(\lambda-\omega)t}\left(1 \vee t^{-\gamma}\right)\,dt.
\end{equation}
and by choosing $\lambda_{1}>\lambda_{0}\vee \omega$ we get that $F$ is a contraction on $B_{1}$. In particular, this imply that $v\in C^{1}_{b}(\overline{H})$ and by using Lemma (\ref{lemma-mild-fixed}) that also $\nabla^{B}v\in C^{1}_{b}(\overline{H};K)$. Moreover it is easy to check that $\nabla v$ belongs to $S^{P}_{\infty}(\overline{H};\overline{H}')$ and so by an obvious generalization of Proposition (\ref{prop-lift-partial-smoothing}) we deduce that both $\nabla \nabla ^{B}v$ and $\nabla^{B}\nabla v$ are well defined and they must coincide.
\qed

\section{Verification and Synthesis of Optimal Control}\label{sec:verification_and_synthesis}
\subsection{The Verification Theorem}\label{sec:verification_theorem}
This section represents the culmination of our theoretical development. Having established the existence, uniqueness, and regularity of a mild solution $v$ to the HJB equation, we now connect this analytical object back to the original stochastic control problem. The main goal is to prove a "verification theorem". This type of theorem provides a sufficient condition for optimality and, crucially, identifies the solution of the HJB equation with the value function of the control problem, i.e., $v(x) = V(x)$. For other works on verification theorems, see e.g. \cite{StannatWessels}.
\\
To achieve this, we first need to show that our mild solution can be approximated by a sequence of more regular, "classical" solutions. This is a common technique in the theory of viscosity solutions and for mild solutions of HJB equations in infinite dimensions, as direct application of tools like Itô's formula to the mild solution $v$ is not possible due to its limited regularity \cite[Section 4.5]{FabbriGozziSwiech}, \cite[Section 4]{FGFM-II}, \cite[Section 5.1]{GozziMasiero23}. 
\\
Once the approximation result is in place, we will prove the "fundamental identity" (Proposition \ref{prop-fund-id}), a key relationship that holds for any admissible control. This identity will be the foundation of the verification theorem (Theorem \ref{teorema-verifica-suff}), which ultimately allows us to confirm the optimality of a given control strategy.

\begin{definition}\label{def-diff-operator}
We define the following subset of $UC_{b}(\overline{H})$
\begin{equation*}
\cald(\cala_{0})= \left\{\phi\in UC^{2}_{b}(\overline{H}):\quad \bar A^{*}\nabla\phi\in UC_{b}(\overline{H};\overline{H}'), \quad Q\nabla ^2\phi\in UC_{b}(\overline{H}; \call_{1}(\overline{H})) \right\}
\end{equation*}
and the operator $\cala_{0}:\cald(\cala_{0})\rightarrow UC_{b}(\overline{H})$ as
\begin{equation*}
\cala_{0}[\phi](x)= \frac{1}{2} \operatorname{Tr} \left(Q\nabla ^2\phi(x)\right)
+ \<\bar A^*\nabla \phi(x),x\>_{\overline{H}',\overline{H}}.
\end{equation*}
\end{definition}
We defined the domain of the infinitesimal generator $\mathcal{A}_0$ in a rigorous way. Note that the condition $\nabla\phi \in D(\bar{A}^*)$ is a strong regularity requirement, typical for defining classical solutions of HJB equations in infinite dimensions {\cite[Section 4]{DaPratoZabczyk02}}, {\cite{Cerrai95}}.

\begin{definition}\label{def-classical-sol}
Let $g\in C_{b}(\overline{H})$. A function $v:\overline{H}\rightarrow \R$ is a classical solution of
\begin{equation} \label{eq-hjb-g}
\lambda v(x)-\cala_{0}[v](x)=\ell_0(x)+
H_{min}\left(\nabla^B v(x)\right)+g(x)
\end{equation}
if $v\in \cald(\cala_{0})$, $\nabla^B v\in C_{b}(\overline{H};K)$ and satisfies equation (\ref{eq-hjb-g}) pointwise for all $x\in \overline{H}$.
\end{definition}

\begin{definition}\label{def-k-conv} Let $\overline{H}$ be a real and separable Banach space.
A sequence $(f_n)_{n\geq 0}\in C_b(\overline{H})$
is said to be $\calk$-convergent to a function $f\in C_b(\overline{H})$ (and we shall write
$f_n\overset{\calk}{\rightarrow} f$ or $f=\calk-\lim_{n\rightarrow\infty}f_n$) if for any compact set
$\calk\subset \overline{H}$
\vspace{-0.3truecm}
\[
\sup_{n\in\N}\left\| f_n\right\|_{C_b(\overline{H})}<+\infty \quad {\rm and } \quad
\lim_{n\rightarrow\infty}\sup_{x\in\calk}\left| f(x)-f_n(x)\right| =0.
\]
\end{definition}
This is the notion of convergence on compact sets, which is weaker than uniform convergence but strong enough for many purposes in the analysis of SPDEs. It is particularly useful when dealing with semigroups that are not uniformly continuous but are continuous on compact sets {\cite{Cerrai95}}.

\begin{definition}\label{def-strong-sol}
We say that $v:\overline{H}\to\R$ is a $\calk$-strong solution of equation (\ref{HJB-ellittica}) if there exist a sequence $(v_{n})\subset \cald(\cala_{0})$ and $(g_{n})\subset C_{b}(\overline{H})$ such that $v_{n}$ is a classical solution of
\vspace{-0.2truecm}
\begin{equation*}
\lambda w(x)-\cala_{0}[w](x)=\ell_0(x)+
H_{min}\left(\nabla^B w(x)\right)+g_{n}(x)
\vspace{-0.2truecm}
\end{equation*}
and moreover, as $n\to\infty$ the following convergences
\vspace{-0.2truecm}
\begin{equation}
\left\lbrace
\begin{array}{l}
\calk-\lim_{n\rightarrow\infty}v_{n}=v,
\\
\calk-\lim_{n\rightarrow\infty}\nabla^{B}v_{n}= \nabla^{B}v,
\\
\calk-\lim_{n\rightarrow\infty}g_{n}=0.
\end{array}
\right.
\vspace{-0.2truecm}
\end{equation}
\end{definition}
This definition bridges the gap between the mild solution we found and the classical solutions needed for applying Itô's formula. A $\calk$-strong solution is a function that, while not necessarily a classical solution itself, can be approximated by a sequence of classical solutions in a meaningful way \cite[Lemma 4.3]{FGFM-II}.

\begin{proposition}\label{prop-mild-are-strong}
Let hypothesis (\ref{ipo-data}), (\ref{ipo-traj}). (\ref{ipo-growth-paths}) and (\ref{ipo-smoothing-lifting}) be satisfied. Assume moreover either that the extended semigroup $\{\overline{e^{tA}}\}_{t\geq0}$ is analytic or that the image of $P^{*}$ is contained in $\cald({A^{*}})$. Let $v:\overline{H}\to\R$ be the unique mild solution of equation (\ref{HJB-ellittica}) in the class $S^{P}_{\infty}(\overline{H})$. Then, $v$ is also a $\calk$-strong solution in the sense of definition (\ref{def-strong-sol}).
\end{proposition}
\dim
Let $g:\overline{H}\to\R$ be the function defined by
\vspace{-0.2truecm}
\begin{equation*}
g(x):=\ell_{0}(x)+H_{min}\left(\nabla^{B}v(x)\right), \quad x\in \overline{H}.
\vspace{-0.2truecm}
\end{equation*}
Clearly $g$ belongs to $S^{P}_{\infty}(\overline{H})$, hence there exists a measurable function $\hat g: L^{2}_{\rho}\to\R$ such that $g=\hat g\circ \Upsilon^P_\infty$. Following the ideas of \cite[Proof of Lemma 5.4]{GozziMasiero23} we can construct a sequence $(\hat g_{n})_{n\in\N}\subset \calf C^{\infty}_{0}(L^{2}_{\rho})$ such that $\calk-\lim_{n\to\infty} \hat g_{n}=\hat g$. We fix $g_{n}:= \hat g_{n}\circ\Upsilon^P_\infty$ and we still get the convergence $\calk-\lim_{n\to\infty} g_{n}= g$. Now, we define the sequence of aproximate solutions as
\vspace{-0.2truecm}
\begin{equation}
v_{n}(x)= \int_{0}^{\infty} e^{-\lambda t} P_{t}[g_{n}(\cdot)](x)\,dt, \quad x\in \overline{H}, n\in \N.
\vspace{-0.2truecm}
\end{equation}
By definition, $v_{n}$ is set to be the mild solution of the following elliptic equation on $\overline{H}$
\vspace{-0.2truecm}
\begin{equation}\label{appmild}
\lambda w(x)-\cala_{0}[w](x)= g_{n}(x).
\vspace{-0.2truecm}
\end{equation}
We want to show that $v_{n}$ is actually a classical solution of equation (\ref{appmild}). Since $\hat g_{n}\in\calf C^{\infty}_{0}(L^{2}_{\rho})$ it is immediate to verify that $g_{n}\in UC_{b}^{2}(\overline{H})$ and $Q\nabla^{2}g_{n}\in UC_{b}(\overline{H};\call_{1}(\overline{H}))$. We still need to verify that $\nabla g_{n}(x)\in\cald(\overline{A}^{*})$ in order to show that $g_{n}\in \cald(\cala_{0})$. First of all we observe that
\vspace{-0.2truecm}
\begin{equation*}
\<\nabla g_{n}(x),h\>_{\overline{H}',\overline{H}}=\<\nabla \hat g_{n}\left(\Upsilon^P_{\infty}x\right), \Upsilon^P_{\infty}h\>_{L^{2}_\rho}=\<\left(\Upsilon^P_{\infty}\right)^{*}\nabla \hat g_{n}\left(\Upsilon^P_{\infty}x\right),h\>_{\overline{H}',\overline{H}}, \quad x,h\in \overline{H}
\vspace{-0.2truecm}
\end{equation*}
Hence, we have the equivalence $\nabla g_{n}(x)=\left(\Upsilon^P_{\infty}\right)^{*}\nabla \hat g_{n}\left(\Upsilon^P_{\infty}x\right)$, for any $x\in \overline{H}$. By using proposition (\ref{prop-adjoint}) we can write
\vspace{-0.2truecm}
\begin{equation}
\nabla g_{n}(x)=\int_{0}^{+\infty}e^{-\rho s}e^{sA^*}P^*\nabla \hat g_{n}\left(\Upsilon^P_{\infty}x\right)(s)\,ds, \quad x\in \overline{H}.
\vspace{-0.2truecm}
\end{equation}
By hypothesis either the semigroup $\overline{e^{tA}}$ is analytic, or the image of $P^{*}$ is contained in $\cald({A^{*}})$. In both cases we can deduce that $\int_{0}^{+\infty}e^{-\rho s}e^{sA^*}P^*\nabla \hat g_{n}\left(\Upsilon^P_{\infty}x\right)(s)\,ds\in\cald(\overline{A^{*}})$ and then $\nabla g_{n}(x)\in\cald(\overline{A^{*}})$ for any $x\in \overline{H}$. Since $v_{n}$ is a mild solution of equation (\ref{appmild}) and belongs to $\cald(\cala_{0})$ we deduce that $v_{n}$ is a classical solution of equation (\ref{appmild}). Now, we prove the convergences as in (\ref{def-k-conv}). By construction we have that $\calk-\lim_{n\to\infty} g_{n}= g$. Moreover, using the $\calk$-convergence of $g_{n}$ to $g$ and the dominated convergence theorem it is easy to check that
\vspace{-0.2truecm}
\begin{equation}
\calk-\lim_{n\to\infty} \int_{0}^{\infty} e^{-\lambda t} P_{t}[g_{n}](\cdot)\,dt= \int_{0}^{\infty} e^{-\lambda t} P_{t}[g](\cdot)\,dt,
\vspace{-0.2truecm}
\end{equation}
so that $\calk-\lim_{n\to\infty}v_{n}=v$. It remains to prove that $\calk-\lim_{n\to\infty}\nabla^{B}v_{n}=\nabla^{B}v$. We observe that $g_{n}\in S^{P}_{\infty}(\overline{H})$ and the same holds for $v_{n}$. Hence, there exists a measurable function $\hat v_{n}: L^{2}_{\rho}\to\R$ such that $v_{n}=\hat v_{n}\circ \Upsilon^{P}_{\infty}$. By deriving under the integral sign we get
\vspace{-0.2truecm}
\begin{equation}
\nabla^{B}v_{n}(x)= \int_{0}^{\infty} e^{-\lambda t} \nabla^{B} P_{t}[g_{n}(\cdot)](x)\,dt, \quad x\in \overline{H}.
\vspace{-0.2truecm}
\end{equation}
By using proposition (\ref{prop-lift-partial-smoothing}) we can express the $B$-differential of $v_{n}$ as
\vspace{-0.2truecm}
\begin{equation}\label{eqbdiffn}
\int_{0}^{\infty} e^{-\lambda t} \int_{L^2_\rho}\hat g_{n}\left(z_1+ \Upsilon^P_\infty x \right)
\<\widehat\Lambda^{P,B}(t) k, (\Upsilon^P_\infty Q_t(\Upsilon^P_\infty)^*)^{-1/2} z_1\>_{L^2_\rho}
\,\caln(0,(\Upsilon^P_\infty Q_t(\Upsilon^P_\infty)^*))(dz_1)\,dt.
\vspace{-0.2truecm}
\end{equation}
In the same way, being $v\in S^{P}_{\infty}(\overline{H})$ we can express the $B$-differential of $v$ as
\vspace{-0.2truecm}
\begin{equation}\label{eqbdiff}
\int_{0}^{\infty} e^{-\lambda t} \int_{L^2_\rho}\hat g\left(z_1+ \Upsilon^P_\infty x \right)
\<\widehat\Lambda^{P,B}(t) k, (\Upsilon^P_\infty Q_t(\Upsilon^P_\infty)^*)^{-1/2} z_1\>_{L^2_\rho}
\,\caln(0,(\Upsilon^P_\infty Q_t(\Upsilon^P_\infty)^*))(dz_1)\,dt.
\vspace{-0.2truecm}
\end{equation}
By using the $\calk$-convergence of $\hat g_{n}$ to $\hat g$ and the dominated convergence theorem, we get that (\ref{eqbdiffn}) $\calk$-converges to (\ref{eqbdiff}).
\qed

\begin{proposition}\label{prop-fund-id}
Let the hypothesis of proposition (\ref{prop-mild-are-strong}) be satisfied. Let $v$ be the unique mild solution of equation (\ref{HJB-ellittica}) in the class $S^{P}_{\infty}(\overline{H})$.
Then for every $x\in \overline{H}$, and for every admissible control $u\in\cal U$, we have the fundamental identity
\vspace{-0.2truecm}
\begin{equation}\label{relfond}
v(x)
=J(x;u)+\E\int_0^\infty e^{-\lambda s}\left[H_{min}\left(\nabla^B v(X(s))\right)
- H_{CV}\left(\nabla^B v(X(s));u(s)\right)\right]\,ds,
\vspace{-0.2truecm}
\end{equation}
where $X(\cdot)$ is the mild solution of equation $(\ref{equazione stato})$ with initial datum $x\in \overline{H}$ and control $u\in\calu$.
\end{proposition}

\dim
The proof follows the lines of the proof of \cite[Proposition 5.1]{FGFM-II}.
In this proof, to avoid heavy notation, we write
$A$ for $\overline{A}$.
\\
Take any admissible state-control couple $(X(\cdot),u(\cdot))$,
and let $v_n: \overline{H}\to \R$ be the approximating sequence
of strict solutions defined in Proposition (\ref{prop-mild-are-strong}). To apply the Ito formula to
$v_n$, we need to regularize $X(\cdot)$ since it does not live in $D(A)$; to this aim we define for $k \in \N$, sufficiently large,
%\vspace{-0.2truecm}
%\begin{equation*}\label{eq:Xkdef}
$X_k(s)=k(k-A)^{-1}X(s)$.
%\vspace{-0.2truecm}
%\end{equation*}
The process $X_k$ is in $D(A)$, it converges to $X$ ($\P$-a.s. and $s\in [0,\infty)$ a.e.)
and it is a strong solution\footnote{Here we mean strong in the probabilistic sense
and also in the sense of {\cite[Chapter 5, Section 5.1.1]{DaPratoZabczyk14}}.}
of the Cauchy problem
\vspace{-0.2truecm}
\begin{equation*}
\left\{
\begin{array}
[c]{l}
dX_k(s)  =AX_k(s)\,ds+B_ku(s)\,ds+G_kdW_s
,\text{ \ \ \ }s\in [0,\infty) \\
X_k(0)= x_k,
\end{array}
\right.
\vspace{-0.2truecm}
\end{equation*}
where $B_k=k(k-A)^{-1}B$, $G_k=k(k-A)^{-1}G$ and $ x_k=k(k-A)^{-1}x$.
Now observe that the operator $B_k$ is continuous in $\overline{H}$, hence we can apply Dynkin's formulam (see \cite[Section 1.7]{FabbriGozziSwiech}) on the interval $[0,T]$ to the process $e^{\lambda s}v_{n}(X_k(s))$, obtaining
\vspace{-0.2truecm}
\begin{align}\label{quasirelfondv^nk}
& \E v_n(X_k(T)) - v_n(x_k)=\E\int_0^{T} e^{-\lambda s}\left(-\lambda v_{n}(X_{k}(s))+ \frac{1}{2}\operatorname{Tr}\left[Q D^{2}v_{n}(X_{k}(s))\right]  \right)\,ds\nonumber\\
& +\E\int_0^{T} e^{-\lambda s}\left(  \<X_{k}(s), A^{*}Dv_{n}(X_{k}(s))\> + \<B_{k}u(s), Dv_{n}(X_{k}(s))\>    \right)\,ds.
\end{align}
Now we pass to the limit for $k\rightarrow \infty$. Applying the Dominated Convergence Theorem to all terms but the last we get
\vspace{-0.2truecm}
\begin{align}\label{quasirelfondv^n}
& \E v_n(X(T)) - v_n(x)=\E\int_0^{T} e^{-\lambda s}\left(-\lambda v_{n}(X(s))+ \frac{1}{2}\operatorname{Tr}\left[Q D^{2}v_{n}(X(s))\right]  \right)\,ds\nonumber\\
& +\E\int_0^{T} e^{-\lambda s}  \<X(s), A^{*}Dv_{n}(X(s))\> \,ds+ \lim_{k\to\infty}\E\int_0^{T}e^{-\lambda s}\<B_{k}u(s), Dv_{n}(X_{k}(s))\>\,ds.
\end{align}
Concerning the last term, again following the lines of \cite[Proof of Proposition 5.1]{FGFM-II}
\vspace{-0.2truecm}
\begin{equation}
\label{eq:limnablaC}
\lim_{k \to + \infty}\E\int_0^T e^{-\lambda s}\<B_k u(s),\nabla v_n(X_k(s))\>\,ds
=
\E\int_0^T\<u(s),\nabla^B v_n(X(s))\>\,ds.
\vspace{-0.2truecm}
\end{equation}
Using the fact that $v_{n}$ is a classical solution of equation (\ref{appmild}) we get
\begin{equation}
\E v_n(X(T)) - v_n(x)=\E\int_0^{T} e^{-\lambda s} \left[g_{n}(X(s))+\<u(s),\nabla^{B} v_n(X(s))\>\right]\,ds.
\end{equation}
Now we let $n\rightarrow\infty$. By Proposition (\ref{prop-mild-are-strong}), we know that
$v_n\rightarrow v$, $\nabla^Bv_n\rightarrow \nabla^B v$ and $g_{n}\to H_{min}(\nabla^{B}v)$
in the sense of $\calk$-convergence, obtaining by dominated convergence theorem
\begin{equation}
\E v(X(T)) - v(x)=\E\int_0^{T} e^{-\lambda s} \left[H_{min}\left(\nabla^{B}v(X(s))\right)+\ell_{0}(X(s))+\<u(s),\nabla^{B} v(X(s))\>\right]\,ds.
\end{equation}
Now, adding and subtracting $\E\dis\int_t^T  \ell_1(u(s))\,ds$ and rearranging the terms, we get
\vspace{-0.2truecm}
\begin{align*}
v(x)&
=\E v(X(T))+\E\int_0^T \left[\ell_0(X(s))+ \ell_1(u(s))\right]\,ds\\
&+\E\int_0^T \left[H_{min}\left(\nabla^B v(X(s))\right)-
H_{CV}\left(\nabla^B v(X(s));u(s)\right)\right]
\,ds
\vspace{-0.2truecm}
\end{align*}
which immediately gives the claim by passing to the limit for $T\to\infty$.
\qed

\begin{theorem}
\label{teorema-verifica-suff}
Let the hypotheses of proposition (\ref{prop-mild-are-strong}) hold true.
Let $v$ be the unique mild solution of equation (\ref{HJB-ellittica}) in the class $S^{P}_{\infty}(\overline{H})$. Then the following holds.
\vspace{-0.2truecm}
\begin{itemize}
\item[(i)] For all $x\in \overline{H}$ we have
$v(x) \le V(x)$, where $V$ is the value function
defined in \eqref{funzione valore}.
\vspace{-0.2truecm}
\item [(ii)]
Let $x\in \overline{H}$ be fixed.
If, for an admissible control $u^*\in \calu$, we
have, calling $X^*$ the corresponding state,
\vspace{-0.2truecm}
\[
u^*(s)\in \arg\min_{u\in U}H_{CV}\left(\nabla^Bv(X^*(s));u\right)
\vspace{-0.2truecm}
\]
for a.e. $s\in [0,\infty)$, $\P$-a.s., then the pair $(u^*,X^*)$ is
optimal for the control problem starting from $x$
and $v(x)=V(x)=J(x;u^*)$.
\end{itemize}
\end{theorem}
\dim
Let $(X(\cdot),u(\cdot))$ be any admissible state-control couple. Since $H_{min}\left(\nabla^B v(X(s))\right)
\leq H_{CV}\left(\nabla^B v(X(s));u(s)\right)$ for all $s\geq0$ it follows by Proposition (\ref{prop-fund-id}) that $v(x)\leq J(x;u)$. By taking the infimum over all admissible controls $u\in \calu$ we get immediately $(i)$. If now $(X^{*}(\cdot),u^{*}(\cdot))$ satisfies $u^*(s)\in \arg\min_{u\in U}H_{CV}\left(\nabla^Bv(X^*(s));u\right)$ then clearly $H_{min}\left(\nabla^B v(X^{*}(s))\right)
=H_{CV}\left(\nabla^B v(X^{*}(s));u^{*}(s)\right)$ for all $s\geq0$. Again by Proposition (\ref{prop-fund-id}) we get that $v(x)=J(x;u^{*})$, which proves $(ii)$.
\qed

\subsection{Synthesis of the Optimal Feedback Control}\label{sec:feedback_control}
%Having established that the solution to the HJB equation is indeed the value function of our problem ($v=V$), 
We now 
%address the final and most practical question: how to construct an optimal control? This section leverages 
prove the verification theorem and the synthesis of optimal controls. %This means the optimal control action at any time $s$ will be a direct function of the current state of the system, $X(s)$.
%\\
%We will first define the feedback map based on the minimization of the Hamiltonian. Then, we will formulate the resulting closed-loop equation. The main challenge is to ensure that this equation is well-posed, i.e., that it has a unique solution. This typically requires some form of Lipschitz continuity on the feedback function, which translates back to regularity assumptions on the value function $v$ and the problem data. This section provides the final piece of the puzzle, delivering a complete solution to the control problem.

\begin{definition}\label{defdiPsi}
We define, for $x\in \overline{H}$,
the {\em feedback map}
\vspace{-0.2truecm}
\begin{equation}
\Psi(x):=\arg \min_{u\in U} H_{CV}\left(\nabla^Bv(x);u\right),
\vspace{-0.2truecm}
\end{equation}
where, as usual, $v$ is the unique mild solution of equation (\ref{HJB-ellittica}) in the class $S^{P}_{\infty}(\overline{H})$.
The feedback map $\Psi(x)$ provides the optimal action to take when the system is in state $x$. It is found by minimizing the Hamiltonian, using the gradient of the value function $\nabla^Bv(x)$ as the "costate".
\end{definition}

Given any $x\in \overline{H}$,
the so-called Closed Loop Equation
is written, formally, as
\vspace{-0.2truecm}
\begin{equation}\label{cleinclusion}
\left\{
\begin{array}
[c]{l}
dX(s) \in AX(s)\,ds+B\Psi\left(X(s)\right)\,ds+G\,dW_s
,\text{ \ \ \ }s\in [0,\infty) \\
X(0)= x.
\end{array}
\right.
\vspace{-0.2truecm}
\end{equation}

\begin{corollary}
\label{cr:optimalfeedback}
Let the assumptions of proposition (\ref{prop-mild-are-strong}) be satisfied.
Let $v$ be the mild solution of (\ref{HJB-ellittica}).
Fix $x\in \overline{H}$ and assume that the map $\Psi$
defined in (\ref{defdiPsi}) admits a measurable selection $\psi$
such that the Closed Loop Equation
\vspace{-0.2truecm}
\begin{equation}
\label{eq:CLEselection}
\left \{
\begin{array}{l}
d X(s) = AX(s)\,ds+B\psi\left(X(s)\right)\,ds+G\,dW_s
,\text{ \ \ \ }s\in [0,\infty) \\
X(0)=x.
\end{array}
\right.
\vspace{-0.2truecm}
\end{equation}
has a mild solution $X_\psi(\cdot; x)$ (in the sense of {\cite[Chapter 7, Theorem 7.2]{DaPratoZabczyk14}}).
Define, for $s\geq0$, $u_\psi (s)=\psi(X_\psi(s; x))$.
Then the couple
$(u_\psi(\cdot),X_\psi(\cdot;x))$ is optimal at
$x$ and $v(x)=V(x)=J(x;u^*)$.
If, finally, $\Psi(x)$ is always a singleton and the mild solution
of (\ref{HJB-ellittica}) is unique,
then the optimal control is unique.
\end{corollary}
\dim
By construction, the control $u_\psi$ satisfies $u_\psi(s)\in\arg \min_{u\in U} H_{CV}\left(\nabla^Bv(X_{\psi}(s));u\right)$ for all $s\geq0$. Hence, by Theorem \ref{teorema-verifica-suff} we have that the couple $(u_{\psi}(\cdot), X_{\psi}(\cdot))$ is optimal at $x\in \overline{H}$. For the uniqueness, we observe that if $(X(\cdot),u(\cdot))$ is an optimal couple, then $V(x)=J(x;u)$. By Theorem \ref{teorema-verifica-suff} $(i)$ we have that $V(x)\geq v(x)$ and using Proposition \ref{prop-fund-id} we obtain
\begin{equation}\label{disfund}
J(x;u)=V(x)\geq v(x)= J(x;u)+\E\int_0^\infty e^{-\lambda s}\left[H_{min}\left(\nabla^B v(X(s))\right)
- H_{CV}\left(\nabla^B v(X(s));u(s)\right)\right]\,ds.
\end{equation}
However, we can construct an optimal couple in feedback form, say $(X_{\psi}(\cdot),u_{\psi}(\cdot))$, hence by the first part of the proof we get that $v(x)=V(x)$, which implies by using equation (\ref{disfund}) that
\begin{equation}
\E\int_0^\infty e^{-\lambda s}\left[H_{min}\left(\nabla^B v(X(s))\right)- H_{CV}\left(\nabla^B v(X(s));u(s)\right)\right]\,ds=0.
\end{equation}
Since $H_{min}\left(\nabla^B v(X(s))\right)\leq H_{CV}\left(\nabla^B v(X(s));u(s)\right)$ for all $s\geq0$ we deduce that
\begin{equation}
H_{min}\left(\nabla^B v(X(s))\right)= H_{CV}\left(\nabla^B v(X(s));u(s)\right), \quad ds\otimes\P-a.e.
\end{equation}
which clearly implies that $u(s)\in \arg \min_{u\in U} H_{CV}\left(\nabla^Bv(X(s));u\right)$. As $\Psi(X(s))$ is always a singleton given by $\psi(X(s))$ it must hold that $u(s)=u_{\psi}(s)$ and by the uniqueness of mild solution we obtain that $X_{\psi}(s)=X(s)$, $ds\otimes\P-a.e.$
\qed

\begin{hypothesis}\label{ipo-meas-sel}
The set-valued map $P(p):=\left\{ u\in U: \<p,u\>+\ell_1(u)= H_{min}(p)\right\}$ is always non empty; moreover it
admits a Lipschitz continuous selection $\gamma$.
\begin{remark}
This is a crucial assumption for ensuring the well-posedness of the closed-loop equation. The existence of a selection $\gamma$ is guaranteed under mild conditions by selection theorems (e.g., Aumann's, see {\cite[Measurable Choice Theorem]{Aumann67}}). However, Lipschitz continuity is a much stronger requirement and depends on the specific structure of the cost $\ell_1$ and the control set $U$ (e.g., if $\ell_1$ is strictly convex and smooth).
\end{remark}
\end{hypothesis}

\begin{theorem}\label{teo su controllo feedback}
Let Hypotheses \ref{ipo-data}, \ref{ipo-traj}, \ref{ipo-growth-paths}, \ref{ipo-smoothing-lifting}, and \ref{ipo-meas-sel} hold true.
Let $v$ be the unique mild solution of equation (\ref{HJB-ellittica}) in the class $S^{P}_{\infty}(\overline{H})$.
Fix any $x\in \overline{H}$.
Assume also that the running cost $\ell_{0}$ is Lipschitz continuous.
Then the following closed loop equation
\begin{equation}\label{cle}
\left\{
\begin{array}
[c]{l}
dX(s) =AX(s)\,ds+B\gamma\left(\nabla^Bv(X(s))\right)\,ds+G\,dW_s
,\text{ \ \ \ }s\in[0,\infty) \\
X(0)= x
\end{array}
\right.
\end{equation}
admits a unique mild solution $X_\gamma(\cdot;x)$ (in the sense of {\cite[Chapter 7, Theorem 7.2]{DaPratoZabczyk14}}); setting, for $s\geq0$,
$u_\gamma(s):=\gamma\left(\nabla^{B} v(X_\gamma(s;x))\right)$,
we obtain an optimal control at $x$ which is unique if $P$ is always a singleton. Moreover $v(x)=V(x)$.
\end{theorem}
\dim
We just need to prove that equation (\ref{cle}) admits a unique mild solution defined for all time $s\geq0$. The rest of the theorem is an immediate consequence of Corollary \ref{cr:optimalfeedback}. For this purpose, we first prove that the map $x\to\gamma(\nabla^{B}v(x))$ is bounded and Lipschitz continuos. By construction, $\gamma$ takes value in $U$, so it is clearly bounded. Moreover we observe that $\gamma$ is Lipschitz thanks to hypothesis \ref{ipo-meas-sel}, while $\nabla^{B}v$ is Lipschitz by Proposition \ref{prop-regolarità}. Hence the map $x\to\gamma(\nabla^{B}v(x))$ is Lipschitz too, being the composition of Lipschitz functions. The existence and uniqueness of a mild solution follows immediately by contraction arguments for stochastic evolution equations with Lipschitz coefficients {\cite[Chapter 7, Theorem 7.2]{DaPratoZabczyk14}}.
\qed

\section{Example 1: Controlled Stochastic Wave Equation}
\label{sec:wave_application}

We consider a system describing the vibration of a membrane $\Omega \subset \R^d$ with a fixed boundary, subject an internal control force $f(t,x)$ and stochastic noise. The evolution of the state $u(t,x)$ is governed by the following stochastic partial differential equation (SPDE):
\vspace{-0.2truecm}
\begin{equation}\label{eq:wave_spde}
\begin{cases}
    \partial_{tt} u(t,x) = c^2 \Delta u(t,x) + f(t,x) + \sigma dW(t,x) & \text{in } (0,\infty) \times \Omega \\
    u(t,x) = 0 & \text{on } (0,\infty) \times \partial\Omega \\
    u(0,x) = u_0(x), \quad \partial_t u(0,x) = v_0(x) & \text{in } \Omega
\end{cases}
\vspace{-0.2truecm}
\end{equation}
where $c>0$ is the propagation speed, $f(t,x)$ is the distributed control in the space $K = L^2(\Omega;\C)$ and $W(t,x)$ is a cylindrical Wiener process on $L^2(\Omega;\C)$. The state space is the Hilbert space $H = H^1_0(\Omega;\C) \times L^2(\Omega;\C)$, with the norm $\|(u,v)\|_H^2 = c^2\|\nabla u\|_{L^2}^2 + \|v\|_{L^2}^2$, which corresponds to the physical energy of the system. In this space, the SPDE \eqref{eq:wave_spde} is rewritten as a first-order evolution equation:
\begin{equation}\label{eq:wave_abstract}
dX(t) = (A X(t) + B f(t))dt + G dW(t), \quad X(0) = (u_0, v_0)^T.
\end{equation}
The system operators are defined as follows:
 The dynamics operator $A$ is defined on the domain $D(A) = (H^2(\Omega;\C) \cap H_0^1(\Omega;\C)) \times H_0^1(\Omega;\C)$. and its action is defined as
$A = \begin{pmatrix} 0 & I \\ c^2\Delta & 0 \end{pmatrix}$.
It is well known that the operator $A$ generates a contraction semigroup $e^{tA}$ on $H$ and moreover there exists a sequence of eigenvalues of $A$, say $\mu_{n}\in \C$ and the corresponding sequence of eigenfunctions $\Phi_{n}\in H$ which form a Riesz basis of $H$ \cite[Chapter 3]{LasieckaTriggiani}. The control Operator $B:K\to H$ acts on the system only through the second component and is defined as the identity operator of $L^{2}(\Omega;\C)$. Hence, in this framework the control operator is clearly bounded. The set of admissible controls is \vspace{-0.2truecm}
\begin{equation*}\label{eq:admcontr2-wave}
\calu:=\left\{
u:[0,\infty)\times \Omega \to U \subseteq K, \text{ progressively measurable}
\right\},
\vspace{-0.2truecm}
\end{equation*}
for a suitable bounded and closed subset $U$. Finally, the noise operator $G:\R^{N}\to H$ is given by $G = \begin{pmatrix} 0 \\ \sigma \end{pmatrix},$ where $\sigma: \R^N \to L^2(\Omega;\C)$ is a Hilbert-Schmidt operator for an $N$-dimensional noise.
The goal of the optimal control problem is to minimize the following discounted cost functional
\begin{equation}\label{eq:wave_cost}
    J(f) = \E \left[ \int_0^\infty e^{-\lambda t} \left[ \ell_{0}(X(t)) + \ell_{1}(f(t))\right]dt \right].
\end{equation}
\begin{remark}
We observe that this problem does not satisfy the standard strong Feller condition (Hypothesis \ref{ipo-smoothing}). This hypothesis would require $\operatorname{Im}(e^{tA}) \subseteq \operatorname{Im}(Q_t^{1/2})$ for $t>0$.
In this setting, the state space $H = H^1_0(\Omega;\C) \times L^2(\Omega;\C)$ is infinite-dimensional. The wave semigroup $e^{tA}$ (for $\gamma\ge0$) is surjective onto $H$ for any $t>0$. Thus, $\operatorname{Im}(e^{tA}) = H$. Critically, the operator $\sigma$ is defined as Hilbert-Schmidt, which implies that $\sigma$ is a compact operator.
This implies that $G$ itself is a compact operator from its domain into $H$. Consequently, $GG^*$ is also compact.
The covariance operator $Q_t = \int_0^t e^{sA}GG^*e^{sA^*} ds$ is an integral of compact operators (since $e^{sA}$ is bounded) and is therefore itself a compact operator. This implies that $Q_t^{1/2}$ is also compact.
By a standard result of functional analysis, a compact operator on an infinite-dimensional space cannot be surjective. Thus, $\operatorname{Im}(Q_t^{1/2})$ is a proper subspace of $H$.
This leads to the contradiction $\operatorname{Im}(e^{tA}) = H \not\subseteq \operatorname{Im}(Q_t^{1/2})$, and Hypothesis \ref{ipo-smoothing} fails. This failure is the primary motivation for introducing the partial smoothing technique, which only requires smoothing in specific directions.
\end{remark}
Here we assume that $\ell_{0}$ is of the form $\ell_{0}(x)=\hat \ell_{0}(Px)$, for a certain function $\ell_{0}$ defined on $H$, where $P$ is the spectral projection operator onto the subspace $V_N = \operatorname{span}\{\Phi_1, \dots, \Phi_N\}$. We assume, moreover, that $P$ and $G$ have the same image given by $V_{N}$. The crucial step is to verify that the system satisfies a key smoothing condition required by the theory (see Hypothesis \ref{ipo-partial-smoothing}), which involves the interaction between the dynamics, control, noise, and projection. A fundamental property in this context is the commutativity of the semigroup and the projection operator (see Remark \ref{remark commutatività}).

\begin{lemma}
Let $A$ be as in \eqref{eq:wave_abstract} and $P$ be the spectral projection onto the subspace $V_N$ spanned by the first $N$ eigenfunctions of $A$. Then, for all $t \ge 0$, $P$ and $e^{tA}$ commute:
$P e^{tA} = e^{tA}P$.
\end{lemma}
\dim
The proof relies on spectral decomposition. Any state $x\in H$ can be written as $x = \sum_{n=1}^{\infty} c_n \Phi_n$, where $A\Phi_n = \mu_n\Phi_n$. By definition, the action of the operators is:
\vspace{-0.2truecm}
\[ 
e^{tA}x = \sum_{n=1}^{\infty} c_n e^{\mu_n t} \Phi_n \quad \text{and} \quad P x= \sum_{n=1}^{N} c_n \Phi_n. 
\vspace{-0.2truecm}
\]
Computing the compositions, we get:
\vspace{-0.2truecm}
\[ (P e^{tA})x = P \left( \sum_{n=1}^{\infty} c_n e^{\mu_n t} \Phi_n \right) = \sum_{n=1}^{N} c_n e^{\mu_n t} \Phi_n. 
\vspace{-0.2truecm}
\]
\[ (e^{tA}P)x = e^{tA} \left( \sum_{n=1}^{N} c_n \Phi_n \right) = \sum_{n=1}^{N} c_n e^{tA}\Phi_n = \sum_{n=1}^{N} c_n e^{\mu_n t} \Phi_n. 
\vspace{-0.2truecm}
\]
This computation shows that the two operators are equal.
\qed \newline
The key condition for the required estimate concerns the ability of the noise to excite all modes of interest.
\begin{theorem}
If we assume that the finite dimensional operator $PGG^{*}P^{*}$ acting on $V_{N}$ is positive definite, then the model \eqref{eq:wave_abstract} satisfies the following estimate for all $t>0$:
\vspace{-0.2truecm}
\begin{equation}
    \left\| (PQ_tP^*)^{-1/2} P e^{tA}B \right\| \le C t^{-1/2},
\vspace{-0.2truecm}
\end{equation}
where $Q_t = \int_0^t e^{sA}GG^*e^{sA^*}ds$ is the covariance operator. This verifies the abstract hypothesis required by the theory (cf. \cite{GozziMasiero23}) with an exponent $\gamma = 1/2$.
\end{theorem}
\dim
The proof is based on estimating the two terms composing the operator separately. The estimate of the first term is in the following Lemma.
\begin{lemma}
Assume that the finite dimensional operator $PGG^{*}P^{*}$ acting on $V_{N}$ is positive definite, then the projected covariance operator $M(t) = PQ_tP^*$ is invertible for all $t>0$. Moreover, there exist constants $C_1 > 0$ and $t_0 > 0$ such that
$\|\left(PQ_tP^*\right)^{-1/2}\| \le C_1 t^{-1/2}$
for all $t \in (0, t_0]$.
\end{lemma}
\textbf{Proof of the Lemma.}
Since $K = PGG^*P^*$ is positive definite its first eigenvalue $\lambda_{min}(K)$ is strictly positive. For $t \to 0$, an asymptotic expansion of $M(t)$ yields $M(t) = t K + O(t^2)$. This implies that the smallest eigenvalue of $M(t)$, $\lambda_{\min}(M(t))$ can be estimated by $\lambda_{\min}(M(t))\geq\frac{t}{2}\lambda_{\min}(K)$ for small $t$. The following estimate follows directly:
\vspace{-0.2truecm}
\[ \|\left(PQ_tP^*\right)^{-1/2}\| = \frac{1}{\sqrt{\lambda_{\min}(M(t))}} \le \frac{1}{\sqrt{\frac{t}{2}\lambda_{\min}(K)}} = C_1 t^{-1/2}.
\qquad\qquad\qquad\qquad\qed
\vspace{-0.2truecm}
\]
\textbf{Conclusion of the proof of the Theorem.}
In order to conclude the proof of the Theorem we just osberve that the operator $Pe^{tA}B$ is bounded and, by using the fact that $e^{tA}$ is a contraction semigroup, we can estimate the norm as 
$\|Pe^{tA}B\|\leq \|P\| \|e^{tA}\|\|B\|\leq \|P\| \|B\|$.
\hfill\qed
\\
Hence all required assumptions (particularly Hypothesis \ref{ipo-partial-smoothing}) are satisfied.

\section{Example 2: Stochastic Heat Equation with Boundary Control}\label{sec:application}
We now show how to apply our results to the optimal boundary control of a stochastic heat equation. In particular we show that our hypotheses (particularly the crucial "lifting smoothing" condition of Hypothesis \ref{ipo-smoothing-lifting}), are satisfied. For related work on boundary control see also {\cite{MouraFathy13}}.

We consider, in an open connected set with smooth boundary $\calo\subseteq \R^d$ the stochastic heat equation with control at the boundary
\vspace{-0.2truecm}
\begin{equation}\label{eqdiri}
\left\{
\begin{array}{l}
\dis
\partial_{t} y (s,\xi)
= \Delta y(s,\xi)+\sigma \,dW(s,\xi), \qquad s\in (0,\infty),\;
\xi\in \calo,
\\\dis
y(0,\xi)=x(\xi),\; \xi\in \calo,
\\\dis
y(s,\xi)= u(s,\xi), \qquad s\in (0,\infty),\;
\xi\in \partial\calo,
\end{array}
\right.
\vspace{-0.2truecm}
\end{equation}
where $\Delta$ is the Laplace operator, $W$ is a cylindrical Wiener process on $L^{2}(\calo)$, $x(\cdot)\in L^{2}(\calo)$ and $u(s,\cdot)\in L^{2}(\partial \calo)$ for each $s\in(0,\infty)$.\newline
We formulate this problem as an abstract evolution equation as in (\ref{equazione stato}). We choose $H:=L^{2}(\calo)$ to be the state space, $A$ to be the Laplace operator with domain $D(A):=H^{2}(\calo)\cap H^1_0(\calo)$. The operator $A$ is self-adjoint and diagonal with strictly negative eigenvalues $\{-\lambda_n\}_{n\in \N}$
(recall that $\lambda_n\sim n^{2/d}$ as $n \to +\infty$).
We can endow $H$ with a complete orthonormal basis
$\{e_n\}_{n\in \N}$
of eigenvectors of $A$. Moreover it is well known (see, for instance, {\cite[Theorem 4.3]{Pazy83}}) that $A$ generates an analytic semigroup $e^{tA}$ on $H$.\newline
The control space is given by $K:=L^{2}(\partial \calo)$, while the set of admissible controls is 
\vspace{-0.2truecm}
\begin{equation*}\label{eq:admcontr2-heat}
\calu:=\left\{
u:[0,\infty)\times \Omega \to U \subseteq K, \text{ progressively measurable}
\right\},
\vspace{-0.2truecm}
\end{equation*}
for a suitable bounded and closed subset $U$.\newline
Here the extended state $\overline{H}$ is set to be the dual space of $V:=\cald((-A)^{3/4+\varepsilon})$ for some $\epsilon>0$ small. We introduce the Dirichelet map (see e.g.\cite{LionsMagenes} for a complete treatment of the topic) $D: L^{2}(\partial\calo)\rightarrow \cald((-A)^{1/4-\varepsilon})$ as the unique weak solution of
\[
\left\{
\begin{array}{l}
\Delta f(\xi)=0,
\qquad \xi\in \calo,
\\
\dis
f(\xi)=a(\xi), \qquad \xi\in \partial\calo.
\end{array}
\right.
\]
for any boundary data $a\in L^{2}(\partial\calo)$.\newline
The Dirichlet map $D$ is the standard tool to convert a boundary condition into an action on the domain. The regularity result $D \in \call(L^{2}(\partial\calo), \cald((-A)^{1/4-\varepsilon}))$ is a classical result from the theory of elliptic PDEs \cite[Appendix A and Chapter 3]{LasieckaTriggiani}.
We follow the ideas of \cite[Appendix C]{FabbriGozziSwiech} to define a suitable notion of solution for (\ref{eqdiri}). 
\begin{definition}
We say that $X(\cdot)$ is a mild solution of (\ref{eqdiri}) if for all $s\geq0$
\begin{equation}\label{mildiri}
X(s)=e^{sA}x-A\int_0^s e^{(s-r)A}D u(r)\,dr +\int_0^s e^{(s-r)A}G\,dW(r).
\end{equation}
\end{definition}
We observe that we can extend the semigroup $\{e^{tA}\}_{t\geq0}$ to a semigroup $\{\overline{e^{tA}}\}_{t\geq0}$ on $\overline{H}$. The extension procedure relies on the property that the semigroup leaves the smaller space $V$ (equipped with the graph norm) invariant, i.e., for every $v \in V$ and $t \ge 0$, we have $e^{tA}v \in V$. The extended semigroup $\{\overline{e^{tA}}\}_{t \ge 0}$ on $\overline{H}$ is then defined by duality. For any element $x \in \overline{H} = V'$, its action is defined through the duality pairing with an arbitrary test element $v \in V$ as follows:
\begin{equation}\label{eq:extension_duality}
\<\overline{e^{tA}} x, v\>_{V', V} := \<x, e^{tA} v\>_{V', V}.
\end{equation}
This definition is well-posed because $e^{tA}v$ remains in $V$. This construction provides a true extension, as for any $x \in H$, its action is consistent with the original semigroup. Indeed, for any $v \in V$:
$$
\<\overline{e^{tA}} x, v\>_{V', V} = \<x, e^{tA} v\>_{V', V} = \<x, e^{tA} v\>_H = \<e^{tA} x, v\>_H.
$$
Since $V$ is dense in $H$, this implies $\overline{e^{tA}}x = e^{tA}x$ for all $x \in H$. The extended operator inherits the semigroup property from $e^{tA}$.
In the same way we extend the operator $A^{3/4+\epsilon}$ to an operator $\overline{A^{3/4+\epsilon}}\in\call(H;\overline{H})$.\newline
By defining $\overline{A}:= \overline{A^{3/4+\epsilon}}A^{1/4-\epsilon}$ we can write
\vspace{-0.2truecm}
\begin{align*}
& A\int_0^s e^{(s-r)A}D u(r)\,dr=\int_0^s Ae^{(s-r)A}D u(r)\,dr= \int_0^s \overline{A^{3/4+\epsilon}}e^{(s-r)A}A^{1/4-\epsilon}D u(r)\,dr\\
& = \int_0^s e^{(s-r)A}\overline{A^{3/4+\epsilon}}A^{1/4-\epsilon}D u(r)\,dr= \int_0^s e^{(s-r)A}\overline{A}D u(r)\,dr.
\vspace{-0.2truecm}
\end{align*}
In this framework problem (\ref{eqdiri}) can be reformulated as the following abstract evolution equation on $\overline{H}$
\vspace{-0.2truecm}
\begin{equation}\label{eqdiri mild}
\left\{
\begin{array}{l}
\dis
d X(s)= AX(s)\,ds -\overline{A}Du(s)\,ds +G\,dW(s), \quad s\in(0,\infty)
\\\dis
X(0)=x\in H,
\end{array}
\right.
\vspace{-0.2truecm}
\end{equation}
where the control operator
\begin{equation}\label{controlop}
B:=-\overline{A}D
\end{equation}
is an element of $\call(K;\overline{H})$. Notice that thanks to Definition \ref{mildiri}, the evolution is preserved on $H$, hence $X(t;x,u)\in H$ for all $t\geq0$, $x\in H$ and $u\in \cal U$.\newline
Now we consider the optimal control problem related to the stochastic heat equation with boundary control in its abstract reformulation, also in order to introduce the class of function on which we study the partial smoothing, defined by means of the operator $P$. For any given $t\in [0,T]$ and $x \in H$, the objective is to minimize, over all control strategies in $\calu$, the following finite horizon cost:
\vspace{-0.2truecm}
\begin{equation}\label{costoastratto}
J(x;u)=\E \left[\int_0^\infty e^{-\lambda s} \left[\ell_0(X(s))+\ell_1(u(s))\right]\,ds\right].
\vspace{-0.2truecm}
\end{equation}
Moreover we take $Q=(-A)^{-2\beta}$ for some $\beta\ge 0$ and
$P$ a projection on a finite dimensional subspace contained in $(-A)^{-\alpha}$ for some $\alpha>\beta+ \frac14$.\newline
The covariance operator $Q_t$ is given by
\vspace{-0.2truecm}
\begin{equation}\label{cov-heat}
Q_t=\int_0^t (-A)^{-2\beta}e^{2sA}\,ds=(-A)^{-2\beta-1}(I-e^{2tA}).
\vspace{-0.2truecm}
\end{equation}
Notice that it can be deduced by the strong Feller property of the heat transition semigroup that $\operatorname{Im}e^{tA}\subseteq\operatorname{Im}Q_t^{1/2}$, see e.g. {\cite[Chapter 9, Section 9.4 and Chapter 11, Section 11.2.2]{DaPratoZabczyk14}} for a comprehensive bibliography. Now we estimate
$\left\| Q_t^{-1/2}\overline{e^{tA}}B\right\|$.

\begin{lemma}\label{Lemma Qt}
Let $Q_t$ be defined in (\ref{cov-heat}). For every $\varepsilon\in \left(0,\dfrac{1}{4}\right)$, we get, for some $C_0>0$,
\vspace{-0.2truecm}
\begin{equation}\label{Q_t-norm}
\left\| Q_t^{-1/2}\overline{e^{tA}}B\right\|
\le C_0 t^{-\frac{5}{4}-\beta -\varepsilon}.
\vspace{-0.2truecm}
\end{equation}
\end{lemma}
Now we introduce the operator $P$.
Let $\alpha>0$, let $v_1,..., v_n\in D((-A)^{\alpha})$ be linearly independent, and let $P$ be the projection on the span of $\<v_1,...,v_n\>$, namely
\vspace{-0.2truecm}
\begin{equation}\label{P-n-gen}
P:H\rightarrow H,\quad Px=\sum_{i=1}^n\<x,v_i\>v_i, \quad \forall x \in H.
\vspace{-0.2truecm}
\end{equation}
We set moreover, noticing that $P=P^*$,
\vspace{-0.2truecm}
\begin{equation}\label{barQ_t-heat}
\bar Q_t:= PQ_tP=P (-A)^{-1-\beta}(I-e^{2tA})P.
\vspace{-0.2truecm}
\end{equation}
Notice that
$P_\alpha:=(-A)^{\alpha}P$ is a continuous operator on $H$. Hence
\[
\overline{Pe^{tA}}B= P e^{tA}
(-A)^{\frac34+\eps}((-A)^{\frac14-\eps} D), \qquad
(\overline{Pe^{tA}}B)^*= ((-A)^{\frac14-\eps}D)^*
(-A)^{\frac34+\eps-\alpha} e^{tA}P_\alpha
\]
\[
\<Q_tP^*x,P^* x\>=
\<(I-e^{2tA}) (-A)^{-1-2\alpha-\beta} P_\alpha x,P_\alpha x\>
\]
The aim now is to verify that $\operatorname{Im}\left(\overline{Pe^{tA}}(-AD)\right)\subset \operatorname{Im}\left(\bar Q_t^{1/2}\right)$ and to estimate $\left\| \bar Q_t^{-1/2}Pe^{tA}(-AD)\right\|$.
\begin{lemma}\label{Lemma-barQt}
Let $\bar Q_t$ be defined in (\ref{barQ_t-heat}).
Let $\alpha>\beta + \frac14$. Then, for $\varepsilon\in\left(0,\frac14\right)$,
\vspace{-0.2truecm}
\begin{equation}\label{barQ_t-norm-heat}
\operatorname{Im}\left(\overline{Pe^{tA}}B\right)\subset \operatorname{Im}\left(\bar Q_t^{1/2}\right), \quad \left\| \bar Q_t^{-1/2}\overline{Pe^{tA}}B\right\| \leq \frac{C}{t^{1-\varepsilon}}.
\vspace{-0.2truecm}
\end{equation}
\end{lemma}
\begin{remark}
These technical lemmas provide the core estimates needed to show that the abstract hypotheses of our theory are met. Lemma \ref{Lemma-barQt} is particularly important as it establishes the bound required by the partial smoothing hypothesis (specifically, a version of Hypothesis \ref{ipo-smoothing-lifting}) with an integrable singularity $t^{\varepsilon-1}$. This is the key technical verification that allows the application of our HJB theory to this problem.
\end{remark}

\begin{hypothesis}\label{hp:BCcost}
\begin{itemize}
\item[(i)] $\ell_{0}:H\to\R$ is measurable and is such that for the finite set of linearly independent vectors defined in equation \ref{P-n-gen} $\{v_{1},\dots,v_{n}\}\subset \cald((-A)^{\alpha})$ with $\alpha>\frac{1}{4}+\epsilon$ and a suitable function $\bar\ell_{0}\in C_{b}(\R^{n})$ one has
\vspace{-0.2truecm}
\begin{equation}
\ell_{0}(x)=\bar\ell_{0}\left(\<x,v_{1}\>_{H},\dots,\<x,v_{n}\>_{H}\right).
\vspace{-0.2truecm}
\end{equation}
\item[(ii)]$\ell_{1}:U\to\R$ is measurable and bounded from below.
\end{itemize}
\end{hypothesis}
This assumption on the cost function $\ell_0$ ensures that it belongs to the class of lifted functions $\cals^P_\infty(\overline{H})$ (in fact, to the simpler class $B_b^P(H) \subset \cals^P_\infty(\overline{H})$ as per Proposition \ref{prop-adjoint}). This is necessary to apply the partial smoothing results derived from the lifting method.

Now we verify that Hypothesis \ref{ipo-growth-paths} holds true. We perform the computation in the case in which $n=1$, but for the cases in which $n>1$ the proof is identical. We can write
\begin{equation}
\left|\overline{Pe^{tA}}x\right|= \left|\< \overline{e^{tA}}x,v\>v\right|\leq \left|v\right|_{H} \left|\< A^{-\alpha}\overline{e^{tA}}x,A^{\alpha}v\>\right|\leq\left|v\right|_{H}\left|A^{\alpha}v\right|_{H}\left|A^{-\alpha}\overline{e^{tA}}x\right|_{H}.
\end{equation}
Now, for $\theta= \frac{3}{4}+\epsilon-\alpha$ we get
\vspace{-0.2truecm}
\begin{equation}
\left|A^{-\alpha}\overline{e^{tA}}x\right|_{H}=\left|A^{\theta}\overline{e^{tA}}A^{-\alpha-\theta}x\right|_{H}\leq C t^{-\theta} \left|A^{-\alpha-\theta}x\right|_{H}\leq C t^{-\theta}\left|x\right|_{H},
\vspace{-0.2truecm}
\end{equation}
which provides the desired estimate.
Finally we verify Hypothesis \ref{ipo-smoothing-lifting}.
\begin{proposition}
Let $\{P_{t}\}_{t\geq0}$ be the Ornstein-Uhlenbeck semigroup associated to the uncontrolled equation (\ref{eqdiri mild}). Let $B$ be the boundary control operator defined in (\ref{controlop}) and let $P$ be the projection defined in (\ref{P-n-gen}). Then, Hypothesis \ref{ipo-smoothing-lifting} is satisfied with $\gamma=1-\delta$ for some $\delta\in\left(0,\frac{1}{4}\right)$.
\end{proposition}
\dim In this proof we consider the case of the projection on the space generated by only one element $v\in \cald((-A)^\alpha)$, namely $P:H\rightarrow H,\; Px=\<x,v\>v$, for all $x \in H,\, P=P^*$, the extension to a map as in \eqref{P-n-gen} being straightforward. \newline In order to prove Hypothesis \ref{ipo-smoothing-lifting}, point (i), we will prove \eqref{eq:NCdualhigh}. Indeed for any $t>0$ and for $z\in L_\rho^2([0,+\infty),H)$
\vspace{-0.2truecm}
\begin{align}\label{conti1heat}
&\<\Upsilon^P_\infty Q_{t} (\Upsilon^P_\infty )^*z,z\>_{L_\rho^2([0,+\infty),H)}
=\int_0^{t}\left|\int_0^{+\infty}e^{-\rho\tau}Q^{1/2}e^{(\tau+r)\overline{A}^*}P^*z(\tau) \,d\tau\right|^2\, dr\nonumber\\
&=\int_0^{t}\left|\int_0^{+\infty}e^{-\rho\tau}Q^{1/2}e^{(\tau+r)\overline{A}^*}\<z(\tau),v\>v \,d\tau\right|^2\, dr=\int_0^{t}\left| Q^{1/2}e^{r\overline{A}^*}\<\int_0^{+\infty}e^{-\rho\tau}e^{\tau\overline{A}^*}z(\tau)\,d\tau\, ,v\>v \right|^2\, dr\nonumber\\
& =\left| Q_{t}^{1/2}P \int_0^{+\infty}e^{-\rho\tau}e^{\tau\overline{A}^*}z(\tau)\,d\tau \right|^2\geq C t^{1-\delta} \left| {B}^* e^{t\overline{A}^*}P\int_0^{+\infty}e^{-\rho\tau}e^{\tau\overline{A}^*}z(\tau)\,d\tau\right|^2
\vspace{-0.2truecm}
\end{align}
where in the last step we have applied \eqref{barQ_t-norm-heat} and \eqref{eq:NCdualhigh}.

In order to conclude it suffices to notice that, recalling also Proposition \ref{prop-adjoint}-(iv),\vspace{-0.2truecm}
\[
\left| (B^*e^{t\overline{A}^*})(\Upsilon^P_\infty)^* z\right|_K^2= \left|{B}^* e^{t\overline{A}^*}P^* \int_0^{+\infty}e^{-\rho\tau}e^{\tau\overline{A}^*}z(\tau) \,d\tau\right|_K^2.
\vspace{-0.6truecm}
\]
\qed

The preceding analysis shows that the assumptions of the abstract framework, including the fundamental 'lifting smoothing' Hypothesis \ref{ipo-smoothing-lifting}, are verified. 
%The applicability of our theoretical results to the boundary control problem is thus established.

\end{document}